\documentclass[final,onefignum,onetabnum]{siamltex}

	\usepackage{graphics}
	\usepackage{amssymb, amsmath, amsfonts}
	\usepackage{tikz}
	\usepackage{float}
	\usepackage{mathabx}
	\usepackage{stmaryrd}
	\usepackage{rotating}
	
\pgfdeclarelayer{edgelayer}
\pgfdeclarelayer{nodelayer}
\pgfsetlayers{edgelayer,nodelayer,main}

\tikzstyle{none}=[inner sep=0pt]
\tikzstyle{nfilled}=[circle,draw=black,fill=white,inner sep=0pt,minimum size=2mm]
\tikzstyle{filled}=[circle,draw=black,fill=black,inner sep=0pt,minimum size=2mm]

\title{Fitted and unfitted domain decomposition using penalty free Nitsche method for the Poisson problem with discontinuous material parameters}

\author{Thomas Boiveau\thanks{Department of Mathematics, University College London, Gower Street, London, UK-WC1E  6BT, United Kingdom; ({\tt thomas.boiveau.12@ucl.ac.uk}).}
}

\begin{document}

\maketitle

\begin{abstract}
In this paper, we study the stability of the non symmetric version of the Nitsche's method without penalty for domain decomposition. The Poisson problem is considered as a model problem. The computational domain is divided into two subdomain that can have different material parameters. In the first half of the paper we are interested in nonconforming domain decomposition, each subdomain is meshed independently of each other. In the second half, we study unfitted domain decomposition, the computational domain has only one mesh and we allow the interface to cut elements of the mesh. The fictitious domain method is used to handle this specificity. We prove $H^1$-convergence and $L^2$-convergence of the error in both cases. Some numerical results are provided to corroborate the theoretical study.
\end{abstract}


%

\pagestyle{myheadings}
\thispagestyle{plain}
\markboth{THOMAS BOIVEAU}{DOMAIN DECOMPOSITION, PENALTY-FREE NITSCHE METHOD}

\section{Introduction}
The two main methods that can be used for the weak enforcement of the boundary and/or interface conditions are the Lagrange multipliers method and the Nitsche's method that is a penalty based method. The Nitsche's method that has been introduced in 1971 \cite{Nitsche_1971_a} is known to have a symmetric and a nonsymmetric version \cite{Freund_1995_a, Hughes_2000_a}. In this work we consider a nonsymmetric penalty free Nitsche's method \cite{Burman_2012_a, Boiveau_2015_a}, this method can be seen as a Lagrange multiplier method, where the Lagrange multiplier has been replaced by the boundary fluxes of the discrete elliptic operator. The method does not have any additional degrees of freedom nor penalty parameter.

The Nitsche's method has been applied to nonconforming domain decomposition with its symmetric and nonsymmetric version by Becker et al. \cite{Becker_2003_a} for the Poisson problem. The method has been extended by Burman and Zunino \cite{Burman_2006_b} using a weighted average of the fluxes at the interface for the advection-diffusion-reaction problem. Several difficulties in the analysis can be handled by taking the right choice of weights (see \cite{Ern_2009_a,Barrau_2012_a,Burman_2006_b}). For unfitted domain decomposition \cite{Hansbo_2002_a, Barrau_2012_a} the interface can cut elements of the mesh, we handle this problem using the fictitious domain approach \cite{Girault_1995_a,Girault_1999_a,Angot_2005_b,Haslinger_2009_a,Burman_2010_a, Burman_2012_a}.

In the second section of this paper, we study the nonconforming domain decomposition where each subdomain are meshed independently. In the third section, we extend the results to unfitted domain decomposition using the fictitious domain method, the fourth section shows a few numerical examples to illustrate the theoretical study.
	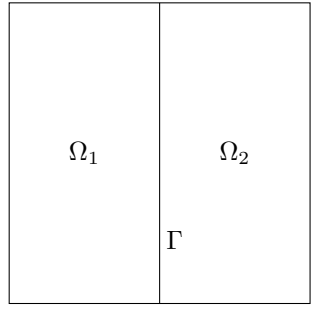
\begin{figure}[h!]
	\begin{center}
	\begin{tikzpicture}[scale=0.5]
	\draw (0,0) rectangle (8cm,8cm);
	\draw[] (4,0)  -- (4,8);
	\draw (4.4,1.7) node {$\Gamma$};
	\draw (2,4) node {$\Omega_1$};
	\draw (6,4) node {$\Omega_2$};
	\end{tikzpicture}
	\end{center}
	\caption{Example of computational domain $\Omega$.}
	\label{domain}
	\end{figure}
	
Let $\Omega_1$ and $\Omega_2$ be two convex bounded domain in $\mathbb{R}^2$ with polygonal boundary, these two domains share an interface $\Gamma=\overline{\Omega}_1\cap\overline{\Omega}_2$. We define the domain $\Omega=\Omega_1\cup\Omega_2$ with boundary $\partial\Omega$, an example of $\Omega$ is represented in Figure \ref{domain}. The Poisson problem considered can be expressed as
\begin{eqnarray}
\label{poisson_DD}
-\mu_i\Delta u^i&=& f_i~~~~\mbox{ in } \Omega_i~,~~i=1,2, \nonumber \\
u^i &=& 0  ~~~~\mbox{ on } \partial\Omega \cap \Omega_i~,~~i=1,2, \nonumber\\
u^1-u^2 &=& 0 ~~~~\mbox{ on } \Gamma,\\
\mu_1\nabla u^1\cdot n_1+\mu_2\nabla u^2\cdot n_2 &=& 0 ~~~~\mbox{ on } \Gamma,\nonumber
\end{eqnarray}
where $u^i$ and $\mu_i$ are respectively the unknown and the diffusivity in $\Omega_i$, $f_i\in L^2\left(\Omega_i\right)$ is a given body force. In this paper $C$ will be used as a generic positive constant that may change at each occurrence, we will use the notation $a\lesssim b$ for $a\leq C b$. For simplicity we will write the $L^2$-norm on a domain $\Theta$, $\left\|\cdot\right\|_{L^2\left(\Theta\right)}$ as $\left\|\cdot\right\|_\Theta$.

\section{Fitted domain decomposition}
\subsection{Preliminaries}
The set $\left\{\mathcal{T}_h^i\right\}_h$ defines the family of quasi-uniform and shape regular triangulations fitted to $\Omega_i$. We define the shape regularity as the existence of a constant $c_\rho\in\mathbb{R}_+$ for the family of triangulations such that, with $\rho_K$ the radius of the largest circle inscribed in an element $K$, there holds
$$\frac{h_K}{\rho_K}\leq c_\rho~~~~\forall K\in \mathcal{T}_h^i, ~~~~i=1,2.$$
In a generic sense we define $K$ as a triangle in a triangulation  $\mathcal{T}_h^i$ and $h_K:=\mbox{diam}(K)$ is the diameter of $K$. Then we define $h_i:=\mbox{max}_{K\in\mathcal{T}_h^i}h_K$ as the mesh parameter for a given triangulation $\mathcal{T}_h^i$. We study the domain decomposition problem with two subdomains that can be meshed independently, we make the assumption $\mu_2h_1\geq\mu_1h_2$ and we set $h:=\text{max}\left(h_1,h_2\right)$.

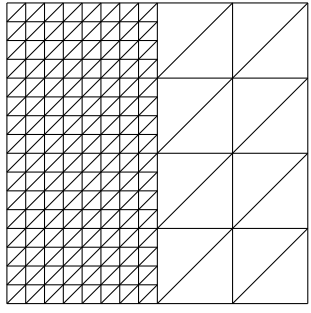
\begin{figure}[h!]
\begin{center}
\begin{tikzpicture}[scale=0.5]
	\begin{pgfonlayer}{nodelayer}
		\node [style=none] (0) at (-4, 4) {};
		\node [style=none] (1) at (-4, -4) {};
		\node [style=none] (2) at (4, -4) {};
		\node [style=none] (3) at (4, 4) {};
		\node [style=none] (4) at (-3.5, 4) {};
		\node [style=none] (5) at (-3, 4) {};
		\node [style=none] (6) at (-2.5, 4) {};
		\node [style=none] (7) at (-2, 4) {};
		\node [style=none] (8) at (-1.5, 4) {};
		\node [style=none] (9) at (-1, 4) {};
		\node [style=none] (10) at (-0.5, 4) {};
		\node [style=none] (11) at (0, 4) {};
		\node [style=none] (12) at (2, 4) {};
		\node [style=none] (13) at (4, 2) {};
		\node [style=none] (14) at (4, 0) {};
		\node [style=none] (15) at (4, -2) {};
		\node [style=none] (16) at (2, -4) {};
		\node [style=none] (17) at (0, -4) {};
		\node [style=none] (18) at (0, -2) {};
		\node [style=none] (19) at (0, 0) {};
		\node [style=none] (20) at (0, 2) {};
		\node [style=none] (21) at (-4, 3.5) {};
		\node [style=none] (22) at (-4, 3) {};
		\node [style=none] (23) at (-4, 2.5) {};
		\node [style=none] (24) at (-4, 2) {};
		\node [style=none] (25) at (-4, 1.5) {};
		\node [style=none] (26) at (-4, 1) {};
		\node [style=none] (27) at (-4, 0.5) {};
		\node [style=none] (28) at (-4, 0) {};
		\node [style=none] (29) at (-4, -0.5) {};
		\node [style=none] (30) at (-4, -1) {};
		\node [style=none] (31) at (-4, -1.5) {};
		\node [style=none] (32) at (-4, -2) {};
		\node [style=none] (33) at (-4, -2.5) {};
		\node [style=none] (34) at (-4, -3) {};
		\node [style=none] (35) at (-4, -3.5) {};
		\node [style=none] (36) at (-3.5, -4) {};
		\node [style=none] (37) at (-3, -4) {};
		\node [style=none] (38) at (-2.5, -4) {};
		\node [style=none] (39) at (-2, -4) {};
		\node [style=none] (40) at (-1.5, -4) {};
		\node [style=none] (41) at (-1, -4) {};
		\node [style=none] (42) at (-0.5, -4) {};
		\node [style=none] (43) at (0, -3.5) {};
		\node [style=none] (44) at (0, -3) {};
		\node [style=none] (45) at (0, -2.5) {};
		\node [style=none] (46) at (0, -1.5) {};
		\node [style=none] (47) at (0, -1) {};
		\node [style=none] (48) at (0, -0.5) {};
		\node [style=none] (49) at (0, 0.5) {};
		\node [style=none] (50) at (0, 1) {};
		\node [style=none] (51) at (0, 1.5) {};
		\node [style=none] (52) at (0, 2.5) {};
		\node [style=none] (53) at (0, 3) {};
		\node [style=none] (54) at (0, 3.5) {};
	\end{pgfonlayer}
	\begin{pgfonlayer}{edgelayer}
		\draw (0.center) to (1.center);
		\draw (1.center) to (2.center);
		\draw (3.center) to (2.center);
		\draw (0.center) to (3.center);
		\draw (11.center) to (17.center);
		\draw (12.center) to (16.center);
		\draw (18.center) to (15.center);
		\draw (19.center) to (14.center);
		\draw (20.center) to (13.center);
		\draw (20.center) to (12.center);
		\draw (19.center) to (3.center);
		\draw (18.center) to (13.center);
		\draw (17.center) to (14.center);
		\draw (16.center) to (15.center);
		\draw (21.center) to (54.center);
		\draw (22.center) to (53.center);
		\draw (23.center) to (52.center);
		\draw (24.center) to (20.center);
		\draw (25.center) to (51.center);
		\draw (26.center) to (50.center);
		\draw (27.center) to (49.center);
		\draw (28.center) to (19.center);
		\draw (29.center) to (48.center);
		\draw (30.center) to (47.center);
		\draw (31.center) to (46.center);
		\draw (32.center) to (18.center);
		\draw (33.center) to (45.center);
		\draw (34.center) to (44.center);
		\draw (35.center) to (43.center);
		\draw (10.center) to (42.center);
		\draw (41.center) to (9.center);
		\draw (8.center) to (40.center);
		\draw (7.center) to (39.center);
		\draw (38.center) to (6.center);
		\draw (5.center) to (37.center);
		\draw (36.center) to (4.center);
		\draw (10.center) to (27.center);
		\draw (26.center) to (9.center);
		\draw (25.center) to (8.center);
		\draw (24.center) to (7.center);
		\draw (23.center) to (6.center);
		\draw (22.center) to (5.center);
		\draw (21.center) to (4.center);
		\draw (28.center) to (11.center);
		\draw (29.center) to (54.center);
		\draw (30.center) to (53.center);
		\draw (31.center) to (52.center);
		\draw (32.center) to (20.center);
		\draw (33.center) to (51.center);
		\draw (34.center) to (50.center);
		\draw (35.center) to (49.center);
		\draw (19.center) to (1.center);
		\draw (36.center) to (48.center);
		\draw (47.center) to (37.center);
		\draw (38.center) to (46.center);
		\draw (18.center) to (39.center);
		\draw (40.center) to (45.center);
		\draw (44.center) to (41.center);
		\draw (42.center) to (43.center);
	\end{pgfonlayer}
\end{tikzpicture}
\caption{Example of mesh of $\Omega$.}
\end{center}
\end{figure}

Let $V_i:=\left\{v\in H^{1}\left(\Omega_i\right):v|_{\partial\Omega}=0\right\}$ for $i=1,2$, then $u=(u^1,u^2)\in V_1\times V_2$.
$\mathbb{P}_k(K)$ defines the space of polynomials of degree less than or equal to $k$ on the element $K$. On each domain $\Omega_i$ we define the space of continuous piecewise polynomial functions
\begin{equation*}
V_{i}^{kh}:=\left\{u_h\in V_i:u_h|_K\in\mathbb{P}_k\left(K\right)~~\forall K\in\mathcal{T}_h^i\right\},~~k\geq 1,
\end{equation*}
$V_h^k=V_{1}^{kh}\times V_{2}^{kh}$, every function in $V_h^k$ has two components $v_h=(v_h^1,v_h^2)$.
We now recall two classical inequalities.
\begin{lemma}
\label{trace}
There exists $C_T\in \mathbb{R}_+$ such that for all $w\in H^1\left(K\right)$ and for all $K\in\mathcal{T}_h$, the \textbf{trace inequality} holds
$$\left\|w\right\|_{\partial K}\leq C_T \left(h_K^{-\frac12}\left\|w\right\|_{K}+h_K^{\frac12}\left\|\nabla w\right\|_{K}\right).$$
\end{lemma}
\begin{lemma}
\label{inverse}
There exists $C_I\in \mathbb{R}_+$ such that for all $u_h\in\mathbb{P}_k(K)$ and for all $K\in\mathcal{T}_h$, the \textbf{inverse inequality} holds
$$\left\|\nabla u_h\right\|_{K}\leq C_I h_K^{-1}\left\|u_h\right\|_{K}.$$
\end{lemma}
At the interface we use the notations
$$\llbracket w \rrbracket=w^1-w^2,$$
for the jump, and
$$\left\{w\right\}=\omega_1w^1+\omega_2w^2~~~,~~~~~~\left<w\right>=\omega_2w^1+\omega_1w^2,$$
with the following weights
$$\omega_1=\frac{h_1\mu_2}{h_1\mu_2+h_2\mu_1}~~~,~~~~\omega_2=\frac{h_2\mu_1}{h_1\mu_2+h_2\mu_1}.$$
At the interface $\Gamma$ we define the normal vector $n:=n_1=-n_2$, then we define
$$\left\{\mu\nabla w \cdot n\right\}=\omega_1\mu_1\nabla w_1 \cdot n+\omega_2\mu_2\nabla w_2 \cdot n,$$
we note that
$$\left\{\mu\nabla u \cdot n\right\}=\mu_1\nabla u^1 \cdot n_1=-\mu_2\nabla u^2 \cdot n_2.$$
To simplify the notations in the analysis we set
$$\gamma=\frac{\mu_1\mu_2}{h_1\mu_2+h_2\mu_1}.$$
We now introduce a structure of patches that will be used in the upcoming inf-sup analysis similarly as \cite{Burman_2012_b, Boiveau_2015_a}. Let the interface elements be the triangles with either a face or a vertex on the interface $\Gamma$. We regroup the interface elements of $\Omega_i$ in closed disjoint patches $P_j^i$ with boundary $\partial P_j^i$, $j=1,...,N_{P}^i$. $N_P^i$ defines the total number of patches. Note that $P^i:=\cup_jP_j^i$. Let $F_j^i:=\partial P_j^i \cap \Gamma$. For each $F_j^i$ there exists two positive constants $c_1$, $c_2$ such that for all $j$
$$c_1h_i\leq\mbox{meas}(F_j^i)\leq c_2 h_i,$$
with $h_i$ the meshsize of the subdomain considered. Let us focus on the patches $P_j^1$ attached to the domain $\Omega_1$. Each patch $P_j^1$ is associated with a function $\phi_j \in V_1^{1h}$ defined such that for each node $x_i\in\mathcal{T}_h^1$ of the face $F_j^1$ let
\begin{equation}
\phi_j\left(x_i\right)=
\left\{
\begin{array}{rllr}
0&\text{for}& x_i \in \Omega_1 \backslash \mathring{F}_j^1\\
1&\text{for}& x_i \in \mathring{F}_j^1,
\end{array}\right.
\label{phi_gamma}
\end{equation}
with $i=1,\dots, N_n$. $N_n$ is the number of node in the triangulation $\mathcal{T}_h^1$ and $\mathring{F}_j^1$ defines the interior of the face $F_j^1$. We define the function $v_\Gamma^1\in V_1^{kh}$ such that
\begin{equation}
v_\Gamma^1=\alpha\sum_{j=1}^{N_P^1} v_j^1,
\label{defvgamma}
\end{equation}
with
\begin{equation}
\label{defvj}
v_j^1=\nu_j\phi_j~~~,~~~~~~\nu_j\in\mathbb{R}.
\end{equation}
\begin{figure}[H]
\begin{center}
\begin{tikzpicture}[scale=0.5]
\draw [black] (-5,0)node{\Large{$\Omega_1$}};
\draw [black] (3,0)node{\Large{$\Omega_2$}};
\draw [black] (0.5,-4)node{\Large{$\Gamma$}};

	\begin{pgfonlayer}{nodelayer}
		\node [style=none] (0) at (0, 6) {};
		\node [style=nfilled] (1) at (-2, -3) {};
		\node [style=nfilled] (2) at (-2, -1) {};
		\node [style=nfilled] (3) at (0, -3) {};
		\node [style=filled] (4) at (0, -1) {};
		\node [style=filled] (5) at (0, 1) {};
		\node [style=nfilled] (6) at (-2, 1) {};
		\node [style=nfilled] (7) at (-2, 3) {};
		\node [style=filled] (8) at (0, 3) {};
		\node [style=nfilled] (9) at (0, 5) {};
		\node [style=none] (10) at (0, -4.5) {};
	\end{pgfonlayer}
	\begin{pgfonlayer}{edgelayer}
		\draw (1.center) to (4.center);
		\draw (1.center) to (3.center);
		\draw (2.center) to (1.center);
		\draw (2.center) to (6.center);
		\draw (6.center) to (5.center);
		\draw (2.center) to (4.center);
		\draw (2.center) to (5.center);
		\draw (6.center) to (7.center);
		\draw (7.center) to (8.center);
		\draw (6.center) to (8.center);
		\draw (7.center) to (9.center);
		\draw (0.center) to (10.center);
	\end{pgfonlayer}
\end{tikzpicture}
\caption{Example Patch $P_j^1$ on the interface elements, the function $\phi_j$ is equal to 0 in the nonfilled nodes, 1 in the filled nodes.}
\end{center}
\end{figure}
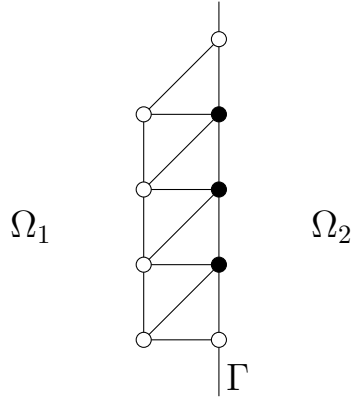
In order to define the properties of $v_j^1$, we define the $P_0$-projection of a function $w$ on an interval $I$
$$\overline{w}^{I}:=\text{meas}\left(I\right)^{-1}\int_{I} w~\text{d}s.$$
We may now define on each face $F_j^1$ the function $v_j^1$ such that
\begin{equation}
\label{vjprop}
\text{meas}\left(F_j^1\right)^{-1}\int_{F_j^1}\nabla v_j^1\cdot n~\text{d}s=h_1^{-1}\overline{\llbracket u_h \rrbracket}^{F_j^1}.
\end{equation}
Applying the Poincar\'e inequality, on each patch $P_j^1$ the function $v_j^1$ has the following property
\begin{equation}
\label{poincare}
\left\|v_j^1\right\|_{P_j^1}\lesssim h_1\left\|\nabla v_j^1\right\|_{P_j^1}.
\end{equation}
It is straightforward to observe that on each face $F_j^i$ it holds
\begin{equation}
\label{stdapprox}
\left\|u_h^i-\overline{u_h^i}^{F_j^i}\right\|_{F_j^i}\lesssim h_i\left\|\nabla u_h^i\right\|_{F_j^i}.
\end{equation}
The Lemma 4.1 of \cite{Burman_2012_b} allows us to write the following inequality for each patch $P_j^1$ of the triangulation $\mathcal{T}_h^1$,
\begin{equation}
\label{keyineq_DD}
\left\|\nabla v_j^1\right\|_{P_j^1}^2\lesssim \left\|h_1^{-\frac12}\overline{\llbracket u_h\rrbracket}^{F_j^1}\right\|_{F_j^1}^2.
\end{equation}
\begin{lemma}
\label{lowbound_projection_DD}
Considering the patches $P_j^i$ as defined above $\forall u_h\in V_h^k$ the following inequality holds
\begin{equation*}
\sum_{j=1}^{N_p^1}\left\|\gamma^\frac12\overline{\llbracket u_h\rrbracket}^{F_j^1}\right\|_{F_j^1}\geq\sum_{j=1}^{N_p^1}\left\|\gamma^\frac12\llbracket u_h\rrbracket\right\|_{F_j^1}-C\omega_1\sum_{j=1}^{N_p^1}\left\|\mu_1^\frac12\nabla u_h^1\right\|_{P_j^1}-C\omega_2\sum_{j=1}^{N_p^2}\left\|\mu_2^\frac12\nabla u_h^2\right\|_{P_j^2}.
\end{equation*}
\end{lemma}
\proof
Considering the triangle inequality and the definition of the jump we can write
\begin{eqnarray*}
\left\|\gamma^\frac12\llbracket u_h\rrbracket\right\|_{F_j^1}
&\leq&\left\|\gamma^\frac12\overline{\llbracket u_h\rrbracket}^{F_j^1}\right\|_{F_j^1}+\left\|\gamma^\frac12\left(\llbracket u_h\rrbracket-\overline{\llbracket u_h\rrbracket}^{F_j^1}\right)\right\|_{F_j^1}\\
&\leq&\left\|\gamma^\frac12\overline{\llbracket u_h\rrbracket}^{F_j^1}\right\|_{F_j^1}+\left\|\gamma^\frac12\left(u_h^1-u_h^2-\left(\overline{ u_h^1}^{F_j^1}-\overline{ u_h^2}^{F_j^1}\right)\right)\right\|_{F_j^1}.
\end{eqnarray*}
Taking the sum over the whole interface and using the triangle inequality once again followed by inequality (\ref{stdapprox}), trace inequality of Lemma \ref{trace} and inverse inequality of Lemma \ref{inverse} we obtain
\begin{eqnarray*}
\sum_{j=1}^{N_p^1}\left\|\gamma^\frac12\llbracket u_h\rrbracket\right\|_{F_j^1}
&\leq&\sum_{j=1}^{N_p^1}\left\|\gamma^\frac12\overline{\llbracket u_h\rrbracket}^{F_j^1}\right\|_{F_j^1}+\gamma^\frac12\sum_{j=1}^{N_p^1}\left\|u_h^1-\overline{ u_h^1}^{F_j^1}\right\|_{F_j^1}+\gamma^\frac12\sum_{j=1}^{N_p^2}\left\|u_h^2-\overline{ u_h^2}^{F_j^2}\right\|_{F_j^2}\\
&\leq&\sum_{j=1}^{N_p^1}\left\|\gamma^\frac12\overline{\llbracket u_h\rrbracket}^{F_j^1}\right\|_{F_j^1}+\gamma^\frac12h_1\sum_{j=1}^{N_p^1}\left\|\nabla u_h^1\right\|_{F_j^1}+\gamma^\frac12h_2\sum_{j=1}^{N_p^2}\left\|\nabla u_h^2\right\|_{F_j^2}\\
&\leq&\sum_{j=1}^{N_p^1}\left\|\gamma^\frac12\overline{\llbracket u_h\rrbracket}^{F_j^1}\right\|_{F_j^1}+C\omega_1\sum_{j=1}^{N_p^1}\left\|\mu_1^\frac12\nabla u_h^1\right\|_{P_j^1}+C\omega_2\sum_{j=1}^{N_p^2}\left\|\mu_2^\frac12\nabla u_h^2\right\|_{P_j^2}.\endproof
\end{eqnarray*}

We end this section by defining the triple norm of a function $w$
\begin{equation}
\label{triplenorm_poisson_DD}
\left\vvvert w \right\vvvert^2=\sum_{i=1}^2\left\|\mu_i^\frac12\nabla w^i\right\|_{\Omega_i}^2+\left\|\gamma^\frac12\llbracket w\rrbracket\right\|_{\Gamma}^2.
\end{equation}

\subsection{Finite element formulation}
Classically for the Poisson problem (\ref{poisson_DD}) for each subdomain domain $\Omega_i$ we obtain by integration by parts
$$\left(\mu_i\nabla u^i,\nabla v^i\right)_{\Omega_i}-\left\langle\mu_i\nabla u^i\cdot n_i,v^i \right\rangle_{\Gamma}=\left(f,v^i\right)_{\Omega_i}~~~,~~~~~~\forall v^i\in V_i.$$
By taking the sum of the boundary terms we obtain
\begin{equation*}
\sum_{i=1}^2\left\langle\mu_i\nabla u^i\cdot n_i,v^i \right\rangle_{\Gamma}
=\int_{\Gamma}\left\llbracket\left(\mu\nabla u \cdot n\right)v\right\rrbracket\text{d}s=\left\langle\left\{\mu\nabla u \cdot n\right\},\left\llbracket v\right\rrbracket\right\rangle_{\Gamma}+\left\langle\left\llbracket\mu\nabla u \cdot n\right\rrbracket,\left< v\right>\right\rangle_{\Gamma},
\end{equation*}
\eqref{poisson_DD} tells us that $\left\llbracket\mu\nabla u \cdot n\right\rrbracket=0$, then we get
$$\sum_{i=1}^2\left(\mu_i\nabla u^i,\nabla v^i\right)_{\Omega_i}-\left\langle\left\{\mu\nabla u \cdot n\right\},\left\llbracket v\right\rrbracket\right\rangle_{\Gamma}=\sum_{i=1}^2\left(f_i,v^i\right)_{\Omega_i}.$$
Adding the Nitsche term, it leads to the following finite element formulation : Find $u_h\in V_h^k$ such that
\begin{equation}
\label{formulationpoisson_DD}
A_h\left(u_h,v_h\right)=L_h\left(v_h\right)~~~\forall v_h\in V_h^k,
\end{equation}
where
\begin{eqnarray*}
A_h\left(u_h,v_h\right)&=&\sum_{i=1}^2\left(\mu_i\nabla u_h^i,\nabla v_h^i\right)_{\Omega_i}-\left\langle\left\{\mu\nabla u_h \cdot n\right\},\left\llbracket v_h\right\rrbracket\right\rangle_{\Gamma}+\left\langle\left\{\mu\nabla v_h \cdot n\right\},\left\llbracket u_h\right\rrbracket\right\rangle_{\Gamma},\\
L_h\left(v_h\right)&=&\sum_{i=1}^2\left(f_i,v_h^i\right)_{\Omega_i}.
\end{eqnarray*}

\subsection{Inf-sup stability}
This section leads to the inf-sup stability of the penalty free scheme previously introduced, we first prove an auxiliary Lemma.
\begin{lemma}
\label{lowerbound_poisson_DD}
For $u_h,v_h\in V_h^k$ with $v_h= u_h+ v_\Gamma^1$, $v_\Gamma^1$ defined by equations \eqref{defvgamma} and \eqref{defvj}, there exists a positive constant $\beta_0$ such that the following inequality holds
$$\beta_0\left\vvvert u_h\right\vvvert^2 \leq A_h(u_h,v_h).$$
\end{lemma}
\proof
Using the definition of $v_\Gamma^1$, we can write the following
$$A_h\left(u_h,v_h\right)=A_h\left(u_h,u_h\right)+ \alpha\sum_{j=1}^{N_p^1}A_h\left(u_h,v_j^1\right).$$
Clearly we have
\begin{equation*}
A_h\left(u_h,u_h\right)=\left\|\mu_1^{\frac12}\nabla u_h^1\right\|_{\Omega_1}^2+\left\|\mu_2^{\frac12}\nabla u_h^2\right\|_{\Omega_2}^2,
\end{equation*}
and
\begin{equation*}
A_h\left(u_h,v_j^1\right)=\left(\mu_1\nabla u_h^1,\nabla v_j^1\right)_{P_j^1}-\left\langle\left\{\mu\nabla u_h\cdot n\right\}, v_j^1\right\rangle_{F_j^1}+\omega_1\left\langle\mu_1 \nabla v_j^1\cdot n,\llbracket u_h\rrbracket\right\rangle_{F_j^1}.
\end{equation*}
Using Cauchy-Schwarz inequality and inequality (\ref{keyineq_DD})
\begin{eqnarray*}
\left(\mu_1\nabla u_h^1,\alpha\nabla v_j^1\right)_{P_j^1}
&\geq&-\left\|\mu_1^\frac12\nabla u_h^1\right\|_{P_j^1}\alpha\mu_1^\frac12\left\|\nabla v_j^1\right\|_{P_j^1}\\
&\geq&-\epsilon\left\|\mu_1^{\frac12}\nabla u_h^1\right\|_{P_j^1}^2-\frac{C\alpha^2}{4\epsilon}\left(1+\frac{h_2\mu_1}{h_1\mu_2}\right)\left\|\gamma^\frac12\overline{\llbracket u_h\rrbracket}^{F_j^1}\right\|_{F_j^1}^2.
\end{eqnarray*}
Using the trace and inverse inequalities, (\ref{poincare}) and (\ref{keyineq_DD}) we can write
\begin{eqnarray*}
&&\left\langle\left\{\mu\nabla u_h\cdot n\right\}, \alpha v_j^1\right\rangle_{F_j^1}=\left\langle\omega_1\mu_1\nabla u_h^1\cdot n+\omega_2\mu_2\nabla u_h^2\cdot n, \alpha v_j^1\right\rangle_{F_j^1}\\
&&=\left\langle\left(\omega_1\mu_1h_1\right)^\frac12\nabla u_h^1\cdot n+\left(\omega_2\mu_2h_2\right)^\frac12\nabla u_h^2\cdot n, \alpha\gamma^\frac12 v_j^1\right\rangle_{F_j^1}\\
&&\leq\left(\left(\omega_1\mu_1h_1\right)^\frac12\left\|\nabla u_h^1\cdot n\right\|_{F_j^1}+\left(\omega_2\mu_2h_2\right)^\frac12\left\|\nabla u_h^2\cdot n\right\|_{F_j^1}\right)\alpha\left\|\gamma^\frac12v_j^1\right\|_{F_j^1}\\
&&\lesssim\left(\left(\omega_1\mu_1h_1\right)^\frac12\left\|\nabla u_h^1\cdot n\right\|_{F_j^1}+\left(\omega_2\mu_2h_2\right)^\frac12\left\|\nabla u_h^2\cdot n\right\|_{F_j^1}\right)\alpha \left\|\gamma^\frac12\overline{\llbracket u_h\rrbracket}^{F_j^1}\right\|_{F_j^1}.
\end{eqnarray*}
Taking the sum over the full boundary $\Gamma$ and using trace and inverse inequalities once again we obtain
\begin{eqnarray*}
\sum_{j=1}^{N_p^1}\left\langle\left\{\mu\nabla u_h\cdot n\right\}, \alpha v_j^1\right\rangle_{F_j^1}
&\leq&\frac{C\alpha^2}{2\epsilon} \sum_{j=1}^{N_p^1}\left\|\gamma^\frac12\overline{\llbracket u_h\rrbracket}^{F_j^1}\right\|_{F_j^1}^2+\epsilon\left(\omega_1\mu_1h_1\right)\sum_{j=1}^{N_p^1}\left\|\nabla u_h^1\cdot n\right\|_{F_j^1}^2\\&&+\epsilon\left(\omega_2\mu_2h_2\right)\sum_{j=1}^{N_p^2}\left\|\nabla u_h^2\cdot n\right\|_{F_j^2}^2\\
&\leq&\frac{C\alpha^2}{2\epsilon} \sum_{j=1}^{N_p^1}\left\|\gamma^\frac12\overline{\llbracket u_h\rrbracket}^{F_j^1}\right\|_{F_j^1}^2+\epsilon\omega_1\sum_{j=1}^{N_p^1}\left\|\mu_1^\frac12\nabla u_h^1\right\|_{P_j^1}^2\\&&+\epsilon\omega_2\sum_{j=1}^{N_p^2}\left\|\mu_2^\frac12\nabla u_h^2\right\|_{P_j^2}^2.
\end{eqnarray*}
Using the property \eqref{vjprop} of $v_j^1$ we can write for each face $F_j^1$
\begin{equation*}
\alpha\omega_1\left\langle\mu_1\nabla v_j^1\cdot n,\llbracket u_h\rrbracket\right\rangle_{F_j^1}
\geq\alpha \left\|\gamma^\frac12\overline{\llbracket u_h\rrbracket}^{F_j^1}\right\|_{F_j^1}^2+\alpha\omega_1\left\langle\mu_1\nabla v_j^1\cdot n,\llbracket u_h\rrbracket-\overline{\llbracket u_h\rrbracket}^{F_j^1}\right\rangle_{F_j^1}.
\end{equation*}
Using the trace inequality and inequality (\ref{keyineq_DD}) we get
\begin{equation*}
\begin{split}
\alpha\omega_1&\left\langle\mu_1\nabla v_j^1\cdot n,\llbracket u_h\rrbracket-\overline{\llbracket u_h\rrbracket}^{F_j^1}\right\rangle_{F_j^1}\\
&\leq \alpha\omega_1\mu_1\left\|\nabla v_j^1\cdot n\right\|_{F_j^1}\left\|\llbracket u_h\rrbracket-\overline{\llbracket u_h\rrbracket}^{F_j^1}\right\|_{F_j^1}\\
&\lesssim\alpha\omega_1\mu_1h_1^{-\frac12}\left\|\nabla v_j^1\right\|_{P_j^1}\left\|u_h^1-u_h^2-\left(\overline{ u_h^1}^{F_j^1}-\overline{ u_h^2}^{F_j^1}\right)\right\|_{F_j^1}\\
&\lesssim \alpha\left\|\gamma^\frac12\overline{\llbracket u_h \rrbracket}^{F_j^1}\right\|_{F_j^1}\left(\frac{\omega_1\mu_1}{h_1}\right)^\frac12\left\|\left(u_h^1-\overline{ u_h^1}^{F_j^1}\right)-\left(u_h^2-\overline{ u_h^2}^{F_j^1}\right)\right\|_{F_j^1}\\
&\lesssim \alpha\left\|\gamma^\frac12\overline{\llbracket u_h \rrbracket}^{F_j^1}\right\|_{F_j^1}\left(\frac{\omega_1\mu_1}{h_1}\right)^\frac12\left(\left\|u_h^1-\overline{ u_h^1}^{F_j^1}\right\|_{F_j^1}+\left\|u_h^2-\overline{ u_h^2}^{F_j^1}\right\|_{F_j^1}\right).
\end{split}
\end{equation*}
Taking the sum over the whole interface $\Gamma$ and using Young's inequality and (\ref{stdapprox})
\begin{eqnarray*}
&&\sum_{j=1}^{N_p^1}\alpha\omega_1\left\langle\mu_1\nabla v_j^1\cdot n,\llbracket u_h\rrbracket-\overline{\llbracket u_h\rrbracket}^{F_j^1}\right\rangle_{F_j^1}\\
&&\leq \sum_{j=1}^{N_p^1}\left[\frac{C\alpha^2}{2\epsilon}\left\|\gamma^\frac12\overline{\llbracket u_h \rrbracket}^{F_j^1}\right\|_{F_j^1}^2+\epsilon\frac{\omega_1\mu_1}{h_1}\left\|u_h^1-\overline{ u_h^1}^{F_j^1}\right\|_{F_j^1}^2+\epsilon\frac{\omega_1\mu_1}{h_1}\left\|u_h^2-\overline{ u_h^2}^{F_j^1}\right\|_{F_j^1}^2\right]\\
&&\leq \frac{C\alpha^2}{2\epsilon}\sum_{j=1}^{N_p^1}\left\|\gamma^\frac12\overline{\llbracket u_h \rrbracket}^{F_j^1}\right\|_{F_j^1}^2+\epsilon\frac{\omega_1\mu_1}{h_1}\sum_{j=1}^{N_p^1}\left\|u_h^1-\overline{ u_h^1}^{F_j^1}\right\|_{F_j^1}^2+\epsilon\frac{\omega_1\mu_1}{h_1}\sum_{j=1}^{N_p^2}\left\|u_h^2-\overline{ u_h^2}^{F_j^2}\right\|_{F_j^2}^2\\
&&\leq \frac{C\alpha^2}{2\epsilon}\sum_{j=1}^{N_p^1}\left\|\gamma^\frac12\overline{\llbracket u_h \rrbracket}^{F_j^1}\right\|_{F_j^1}^2+\epsilon\omega_1\sum_{j=1}^{N_p^1}\left\|\mu_1^\frac12\nabla u_h^1\right\|_{F_j^1}^2+\epsilon\omega_2\sum_{j=1}^{N_p^2}\left\|\mu_2^\frac12\nabla u_h^2\right\|_{F_j^2}^2.
\end{eqnarray*}
It allows us to write
\begin{eqnarray*}
\sum_{j=1}^{N_p^1}\alpha\omega_1\left\langle\mu_1\nabla v_j^1\cdot n,\llbracket u_h\rrbracket\right\rangle_{F_j^1}
&\geq&-\epsilon\omega_1\sum_{j=1}^{N_p^1}\left\|\mu_1^\frac12\nabla u_h^1 \right\|_{P_j^1}^2-\epsilon\omega_2\sum_{j=1}^{N_p^2}\left\|\mu_2^\frac12\nabla u_h^2 \right\|_{P_j^2}^2\\&&+\alpha\left(1-\frac{C\alpha}{2\epsilon}\right) \sum_{j=1}^{N_p^1}\left\|\gamma^\frac12\overline{\llbracket u_h\rrbracket}^{F_j^1}\right\|_{F_j^1}^2.
\end{eqnarray*}
The full bilinear form $A_h$ now has the following lowerbound
\begin{equation*}
\begin{split}
A_h\left(u_h,v_h\right)\geq
&\left\|\mu_1^{\frac12}\nabla u_h^1\right\|_{\Omega_1\backslash P^1}^2+\left\|\mu_2^{\frac12}\nabla u_h^2\right\|_{\Omega_2\backslash P^2}^2+C_a\sum_{j=1}^{N_p^1}\left\|\mu_1^\frac12\nabla u_h^1\right\|_{P_j^1}^2\\
&+C_b\sum_{j=1}^{N_p^2}\left\|\mu_2^\frac12\nabla u_h^2\right\|_{P_j^2}^2+C_c\sum_{j=1}^{N_p^1}\left\|\gamma^\frac12\overline{\llbracket u_h\rrbracket}^{F_j^1}\right\|_{F_j^1}^2,
\end{split}
\end{equation*}
with the constants
\begin{eqnarray*}
C_a&=&1-\epsilon\left(2\omega_1+1\right),\\
C_b&=&1-2\epsilon\omega_2,\\
C_c&=&\alpha\left(1-\alpha \frac{C}{4\epsilon}\left(5+\frac{h_2\mu_1}{h_1\mu_2}\right)\right).
\end{eqnarray*}
Using Lemma \ref{lowbound_projection_DD} it becomes
\begin{equation*}
\begin{split}
A_h\left(u_h,v_h\right)\geq
&\left\|\mu_1^{\frac12}\nabla u_h^1\right\|_{\Omega_1\backslash P^1}^2+\left\|\mu_2^{\frac12}\nabla u_h^2\right\|_{\Omega_2\backslash P^2}^2+\left(C_a-\omega_1CC_c\right)\sum_{j=1}^{N_p^1}\left\|\mu_1^\frac12\nabla u_h^1\right\|_{P_j^1}^2\\
&+\left(C_b-\omega_2CC_c\right)\sum_{j=1}^{N_p^2}\left\|\mu_2^\frac12\nabla u_h^2\right\|_{P_j^2}^2+C_c\sum_{j=1}^{N_p^1}\left\|\gamma^\frac12\llbracket u_h\rrbracket\right\|_{F_j^1}^2.
\end{split}
\end{equation*}
First we fix $\epsilon=\text{min}\left[\frac{1}{2\left(2\omega_1+1\right)},\frac{1}{4\omega_2}\right]$. The constant $C_c$ will be positive for
$$\alpha<\frac{4\epsilon\omega_1}{C\left(4\omega_1+1\right)}.$$
The terms $\left(C_a-\omega_1CC_c\right)$ and $\left(C_b-\omega_2CC_c\right)$ will be both positive for
$$\alpha<\frac{1}{2C}.$$
$\mathcal{O}\left(\epsilon\right)=\text{min}\left[\mathcal{O}\left(1\right),\mathcal{O}\left(\omega_2^{-1}\right)\right]$, $\mathcal{O}\left(\alpha\right)=\text{min}\left[\mathcal{O}\left(\epsilon\omega_1\right),\mathcal{O}\left(1\right)\right]$ and $\mathcal{O}\left(\beta_0\right)=\mathcal{O}\left(\alpha\right)$.\endproof

\begin{theorem}
\label{infsup_poisson_DD}
There exists a positive constant $\beta>0$ such that for all functions $u_h\in V_h^k$ the following inequality holds
\begin{equation*}
\beta\left\vvvert u_h\right\vvvert\leq\underset{v_h\in V_h^k}{\textup{sup}} \frac{A_h\left(u_h,v_h \right)}{\left\vvvert v_h\right\vvvert}.
\end{equation*}
\end{theorem}
\proof
Considering Lemma \ref{lowerbound_poisson_DD}, the only thing that we need to show is
$$\left\vvvert v_h\right\vvvert\lesssim\left\vvvert u_h\right\vvvert.$$
The triangle inequality gives
$$\left\vvvert v_h\right\vvvert=\left\vvvert u_h\right\vvvert+\sum_{j=1}^{N_p^1}\left\vvvert  v_j^1\right\vvvert.$$
The triple norm (\ref{triplenorm_poisson_DD}) gives
\begin{equation*}
\left\vvvert  v_j^1\right\vvvert^2=\left\|\mu_1^\frac12\nabla v_j^1\right\|_{P_j^1}^2+\left\| \gamma^\frac12v_j^1\right\|_{F_j^1}^2.
\end{equation*}
Recalling the inequality (\ref{keyineq_DD}) and
\begin{equation*}
\left\|\gamma^\frac12\overline{\llbracket u_h\rrbracket}^{F_j^1}\right\|_{F_j^1}\lesssim\left\|\gamma^\frac12\llbracket u_h\rrbracket\right\|_{F_j^1}\lesssim \left\vvvert u_h\right\vvvert,
\end{equation*}
it gives the appropriate upper bound
$$\left\|\mu_1^{\frac12} \nabla v_j^1\right\|_{P_j^1}\lesssim\omega_1^{-\frac12} \left\vvvert u_h\right\vvvert.$$
Using the trace inequality of Lemma \ref{trace} and the inequality (\ref{poincare})
$$\left\|\gamma^\frac12v_j^1\right\|_{F_j^1}\lesssim\omega_1^{\frac12}\left\| \mu_1^{\frac12}\nabla v_j^1\right\|_{P_j^1}\lesssim  \left\vvvert u_h\right\vvvert.$$
Note that $\mathcal{O}\left(\beta\right)=\mathcal{O}\left(\beta_0\omega_1^\frac12\right)$.
 \endproof
\subsection{A priori error estimate}
The proof of the stability done in the previous part leads to the study of the a priori error estimate in the triple norm,  the following consistency relation characterizes the Galerkin orthogonality.
\begin{lemma}
\label{galerkin_poisson_DD}
If $u\in H^{2}\left(\Omega_1\right)\times H^{2}\left(\Omega_2\right)$ is the solution of (\ref{poisson_DD}) and $u_h\in V_h^k$ the solution of (\ref{formulationpoisson_DD}) the following property holds
$$A_h\left(u-u_h,v_h\right)=0~~,~~~~\forall v_h\in V_h^k.$$
\end{lemma}
\proof
We observe that $A_h\left(u,v_h\right)=L_h\left(v_h\right)=A_h\left(u_h,v_h\right)$, $\forall v_h\in V_h^k$.
 \endproof

In order to study the a priori error estimate, we introduce an auxiliary norm
$$\left\| w\right\|_\chi= \left\vvvert w\right\vvvert+\left\|\mu_1^\frac12 h_1^{\frac12}\nabla w^1\cdot n \right\|_{\Gamma}+\left\| \mu_2^\frac12 h_2^{\frac12}\nabla w^2\cdot n \right\|_{\Gamma}.$$
\begin{lemma}
\label{triplestar_poisson_DD}
Let $w\in H^{2}\left(\Omega_1\right)\times H^{2}\left(\Omega_2\right)+V_h^k$ and $v_h\in V_h^k$, there exists a positive constant $M$ such that the bilinear form $A_h\left(\cdot,\cdot\right)$ has the property
$$A_h\left(w,v_h\right)\leq M\left\|w\right\|_\chi \left\vvvert v_h\right\vvvert.$$
\end{lemma}
\proof
Using the trace inequality of Lemma \ref{trace} and the Cauchy Schwartz inequality we can show
\begin{eqnarray*}
\left\langle\left\{\mu\nabla w \cdot n\right\},\llbracket v_h\rrbracket\right\rangle_\Gamma
&=&\left\langle\left(\omega_1\mu_1h_1\right)^{\frac12}\nabla w \cdot n+\left(\omega_2\mu_2h_2\right)^{\frac12}\nabla w \cdot n,\gamma^\frac12\llbracket v_h\rrbracket\right\rangle_\Gamma\\
&\leq&\left(\omega_1^\frac12\left\|\mu_1^\frac12h_1^\frac12\nabla w^1\cdot n\right\|_{\Gamma}+\omega_2^\frac12\left\|\mu_2^\frac12h_2^\frac12\nabla w^2\cdot n\right\|_{\Gamma}\right)\left\|\gamma^\frac12\llbracket v_h\rrbracket\right\|_{\Gamma},
\end{eqnarray*}
\begin{eqnarray*}
\left\langle\left\{\mu\nabla v_h \cdot n\right\},\llbracket w\rrbracket\right\rangle_\Gamma
&\leq&\left(\omega_1^\frac12\left\|\mu_1^\frac12h_1^\frac12\nabla v_h^1\cdot n\right\|_{\Gamma}+\omega_2^\frac12\left\|\mu_2^\frac12h_2^\frac12\nabla v_h^2\cdot n\right\|_{\Gamma}\right)\left\|\gamma^\frac12\llbracket w\rrbracket\right\|_{\Gamma}\\
&\lesssim&\left(\omega_1^\frac12\left\|\mu_1^\frac12\nabla v_h^1\right\|_{\Omega_1}+\omega_2^\frac12\left\|\mu_2^\frac12\nabla v_h^2\right\|_{\Omega_2}\right)\left\|\gamma^\frac12\llbracket w\rrbracket\right\|_{\Gamma}.
\end{eqnarray*}
Using these two upper bound it is straightforward to show that
\begin{equation*}
\sum_{i=1}^{2}\left(\mu_i\nabla w^i,\nabla v_h^i\right)_{\Omega_i}+\left\langle\left\{\mu\nabla w \cdot n\right\},\llbracket v_h\rrbracket\right\rangle_\Gamma+\left\langle\left\{\mu\nabla v_h \cdot n\right\},\llbracket w\rrbracket\right\rangle_\Gamma\lesssim\left\| w\right\|_\chi \left\vvvert v_h\right\vvvert.
\end{equation*}
$\mathcal{O}\left(M\right)=\text{min}\left[\mathcal{O}\left(\omega_1^\frac12\right),\mathcal{O}\left(\omega_2^\frac12\right)\right]$.
 \endproof
\begin{proposition}
\label{bounderror_poisson_DD}
If $u\in H^{k+1}\left(\Omega_1\right)\times H^{k+1}\left(\Omega_2\right)$ is the solution of (\ref{poisson_DD}) and $u_h\in V_h^k$ the solution of (\ref{formulationpoisson_DD}), then there holds
$$ \left\vvvert u-u_h\right\vvvert \leq C_\mu h^k\left|u\right|_{H^{k+1}\left(\Omega\right)},$$
with $C_\mu$ a positive constant that depends on $\mu$ and the mesh geometry.
\end{proposition}
\proof
Let $\pi_h^k$ denote the nodal interpolant, the approximation properties for each $K\in\mathcal{T}_h^i$ with $i=1,2$, give
\begin{equation*}
\left\|u-\pi_h^ku\right\|_{K}+h_K\left\|\nabla\left(u-\pi_h^ku\right)\right\|_{K}+h_K^2\left\|D^2\left(u-\pi_h^ku\right)\right\|_{K}\lesssim h_K^{k+1}\left| u \right|_{H^{k+1}(K)}.
\end{equation*}
Using this property and the trace inequality it is straightforward to show that
\begin{equation*}
\left\vvvert u-\pi_h^ku\right\vvvert \lesssim \left\|u-\pi_h^ku\right\|_\chi
\lesssim \mu_1^\frac12h_1^{k}\left|u^1\right|_{H^{k+1}\left(\Omega_1\right)}+\mu_2^\frac12 h_2^{k}\left|u^2\right|_{H^{k+1}\left(\Omega_2\right)}.
\end{equation*}
Using the Galerkin orthogonality of Lemma \ref{galerkin_poisson_DD}, the Theorem \ref{infsup_poisson_DD} and the Lemma \ref{triplestar_poisson_DD} we can write
$$\beta\left\vvvert u_h-\pi_h^ku\right\vvvert\leq\frac{A_h\left(u-\pi_h^ku,v_h\right)}{\left\vvvert v_h\right\vvvert }\leq M\left\|u-\pi_h^ku\right\|_\chi.$$
Using this property and the triangle inequality we obtain
$$\left\vvvert u-u_h\right\vvvert\leq\left\vvvert u-\pi_h^ku\right\vvvert+\frac{M}{\beta} \left\|u-\pi_h^ku\right\|_\chi.$$
The result follows. $\mathcal{O}\left(C_\mu\right)=\mathcal{O}\left(\left(\mu_1^\frac12+\mu_2^\frac12\right)\left(1+\frac{M}{\beta}\right)\right)$.
 \endproof
\begin{lemma}
\label{l2cvg_poisson_DD}
Let $u\in H^{k+1}\left(\Omega_1\right)\times H^{k+1}\left(\Omega_2\right)$ be the solution of (\ref{poisson_DD}) and $u_h$ the solution of (\ref{formulationpoisson_DD}), then there holds
$$\left\|u-u_h\right\|_{\Omega}\leq C_{\mu}' h^{k+\frac12}\left|u\right|_{H^{k+1}\left(\Omega\right)},$$
with $C_{\mu}'$ is a positive constant that depends on $\mu$ and the mesh geometry.
\end{lemma}
\proof
Let $z$ satisfy the adjoint problem
\begin{eqnarray*}
-\mu_i\Delta z^i&=& u^i-u_h^i ~~~~\mbox{ in } \Omega_i~,~~i=1,2, \nonumber \\
z^i &=& 0  ~~~~~~~~~~~\mbox{ on } \partial\Omega \cap \Omega_i~,~~i=1,2, \nonumber\\
z^1-z^2 &=& 0 ~~~~~~~~~~~\mbox{ on } \Gamma,\\
\mu_1\nabla z^1\cdot n_1+\mu_2\nabla z^2\cdot n_2 &=& 0 ~~~~~~~~~~~\mbox{ on } \Gamma,\nonumber
\end{eqnarray*}
Then we can write
\begin{eqnarray*}
\left\|u^i-u_h^i\right\|_{\Omega_i}^2
&=&\left(u^i-u_h^i,-\mu_i\Delta z^i\right)_{\Omega_i}\\
&=&\left(\nabla\left(u^i-u_h^i\right),\mu_i\nabla z^i\right)_{\Omega_i}-\left\langle\mu_i\nabla z^i\cdot n_i, u^i-u_h^i\right\rangle_{\Gamma},
\end{eqnarray*}
it leads to
\begin{eqnarray*}
\left\|u-u_h\right\|_{\Omega}^2
&\leq&\left\|u^1-u_h^1\right\|_{\Omega_1}^2+\left\|u^2-u_h^2\right\|_{\Omega_2}^2,\\
&\leq&\sum_{i=1}^2\left(\nabla\left(u^i-u_h^i\right),\mu_i\nabla z^i\right)-\left\langle\left\{\mu\nabla z\cdot n\right\},\left\llbracket u-u_h\right\rrbracket\right\rangle_{\Gamma},\\
&\leq&A_h\left(u-u_h,z\right)-2\left\langle\left\{\mu\nabla z\cdot n\right\},\left\llbracket u-u_h\right\rrbracket\right\rangle_{\Gamma}.
\end{eqnarray*}
Using the global trace inequality $\left\|\nabla z^i\cdot n\right\|_{\Gamma}\lesssim\left\| z^i\right\|_{H^2\left(\Omega_i\right)}$ for $i=1,2$, we can write
\begin{equation*}
\begin{split}
|\langle\left\{\mu\nabla z\cdot n\right\}&,\left\llbracket u-u_h\right\rrbracket\rangle_{\Gamma}|\\
&\leq\left(\left(\omega_1\mu_1h_1\right)^\frac12\left\|\nabla z^1\cdot n\right\|_{\Gamma}+\left(\omega_2\mu_2h_2\right)^\frac12\left\|\nabla z^2\cdot n\right\|_{\Gamma}\right)\left\|\gamma^\frac12\left\llbracket u-u_h\right\rrbracket\right\|_{\Gamma}\\
&\lesssim\left(\left(\omega_1\mu_1\right)^\frac12+\left(\omega_2\mu_2\right)^\frac12\right)h^\frac12\left\|z\right\|_{H^2\left(\Omega\right)}\left\vvvert u-u_h\right\vvvert.
\end{split}
\end{equation*}
Using Lemma \ref{galerkin_poisson_DD} and $\left(z^i-\pi_h^1z^i\right)|_{\Gamma}\equiv 0$ for $i=1,2$
\begin{eqnarray*}
A_h\left(u-u_h,z\right)
&=&A_h\left(u-u_h,z-\pi_h^1z\right),\\
&=&\sum_{i=1}^2\left(\nabla\left(u^i-u_h^i\right),\mu_i\nabla \left(z^i-\pi_h^1z^i\right)\right)_{\Omega_i}\\
&&+\left\langle\left\{\mu\nabla \left(z-\pi_h^1z\right)\cdot n\right\},\left\llbracket u-u_h\right\rrbracket\right\rangle_\Gamma,\\
&\lesssim&\left((\omega_1\mu_1)^\frac12+(\omega_2\mu_2)^\frac12\right)h\left|z\right|_{H^2\left(\Omega\right)}\left\vvvert u-u_h\right\vvvert.
\end{eqnarray*}
Then we get
$$\left\|u-u_h\right\|_{\Omega}^2\lesssim C_\mu\left(\left(\omega_1\mu_1\right)^\frac12+\left(\omega_2\mu_2\right)^\frac12\right)(h+h^\frac12)h^{k}\left|u\right|_{H^{k+1}\left(\Omega\right)}\left\|z\right\|_{H^2\left(\Omega\right)}.$$
We conclude by applying the regularity estimate $\left\|z\right\|_{H^2\left(\Omega\right)}\lesssim\left\|u-u_h\right\|_{\Omega}$.
$\mathcal{O}\left(C_\mu'\right)=\mathcal{O}\left(C_\mu\left(\left(\omega_1\mu_1\right)^\frac12+\left(\omega_2\mu_2\right)^\frac12\right)\right)$.
 \endproof

\section{Unfitted domain decomposition}
\subsection{Preliminaries}
In this section the two subdomains $\Omega_1$ and $\Omega_2$ are both meshed with only one triangulation. Let $\left\{\mathcal{T}_h\right\}_h$ be a family of quasi-uniform and shape regular triangulation fitted to $\Omega$, the shape regularity is defined in the same way as in the previous section. In a generic sense, $F$ with normal $n_F$ denotes a face of a triangle $K$ in the triangulation $\mathcal{T}_h$. The mesh size is defined as $h:=\mbox{max}_{K\in\mathcal{T}_h}h_K$.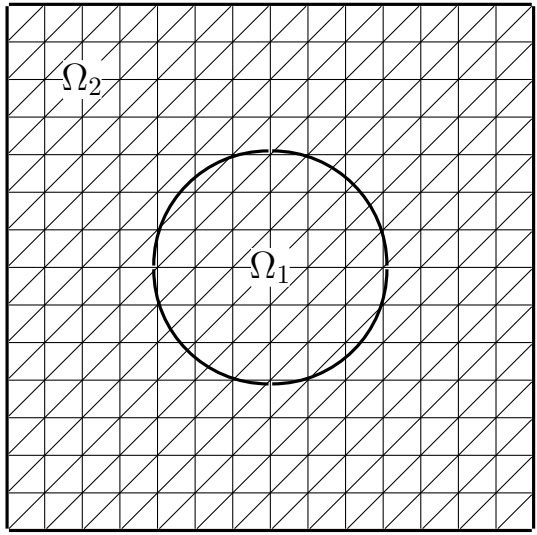
\begin{figure}[H]
\begin{center}
\begin{tikzpicture}[scale=0.5]

\draw [fill=white,white] (-0.5,-0.5) rectangle (0.5,0.5); 
\draw [fill=white,white] (-5.5,4.5) rectangle (-4.5,5.5);

\draw [black] (-5,5)node{\Large{$\Omega_2$}};
\draw [black] (0,0)node{\Large{$\Omega_1$}};

	\begin{pgfonlayer}{nodelayer}
		\node [style=none] (0) at (0, 7) {};
		\node [style=none] (1) at (-1, 7) {};
		\node [style=none] (2) at (-2, 7) {};
		\node [style=none] (3) at (-3, 7) {};
		\node [style=none] (4) at (-4, 7) {};
		\node [style=none] (5) at (-5, 7) {};
		\node [style=none] (6) at (-6, 7) {};
		\node [style=none] (7) at (-7, 7) {};
		\node [style=none] (8) at (1, 7) {};
		\node [style=none] (9) at (2, 7) {};
		\node [style=none] (10) at (3, 7) {};
		\node [style=none] (11) at (4, 7) {};
		\node [style=none] (12) at (5, 7) {};
		\node [style=none] (13) at (6, 7) {};
		\node [style=none] (14) at (7, 7) {};
		\node [style=none] (15) at (-7, 6) {};
		\node [style=none] (16) at (-7, 5) {};
		\node [style=none] (17) at (-7, 4) {};
		\node [style=none] (18) at (-7, 3) {};
		\node [style=none] (19) at (-7, 2) {};
		\node [style=none] (20) at (-7, 1) {};
		\node [style=none] (21) at (-7, 0) {};
		\node [style=none] (22) at (-7, -1) {};
		\node [style=none] (23) at (-7, -2) {};
		\node [style=none] (24) at (-7, -3) {};
		\node [style=none] (25) at (-7, -4) {};
		\node [style=none] (26) at (-7, -5) {};
		\node [style=none] (27) at (-7, -6) {};
		\node [style=none] (28) at (-7, -7) {};
		\node [style=none] (29) at (-6, -7) {};
		\node [style=none] (30) at (-5, -7) {};
		\node [style=none] (31) at (-4, -7) {};
		\node [style=none] (32) at (-3, -7) {};
		\node [style=none] (33) at (-2, -7) {};
		\node [style=none] (34) at (-1, -7) {};
		\node [style=none] (35) at (0, -7) {};
		\node [style=none] (36) at (1, -7) {};
		\node [style=none] (37) at (2, -7) {};
		\node [style=none] (38) at (3, -7) {};
		\node [style=none] (39) at (4, -7) {};
		\node [style=none] (40) at (5, -7) {};
		\node [style=none] (41) at (6, -7) {};
		\node [style=none] (42) at (7, -7) {};
		\node [style=none] (43) at (7, -6) {};
		\node [style=none] (44) at (7, -5) {};
		\node [style=none] (45) at (7, -4) {};
		\node [style=none] (46) at (7, -3) {};
		\node [style=none] (47) at (7, -2) {};
		\node [style=none] (48) at (7, -1) {};
		\node [style=none] (49) at (7, 0) {};
		\node [style=none] (50) at (7, 1) {};
		\node [style=none] (51) at (7, 2) {};
		\node [style=none] (52) at (7, 3) {};
		\node [style=none] (53) at (7, 4) {};
		\node [style=none] (54) at (7, 5) {};
		\node [style=none] (55) at (7, 6) {};
		\node [style=none] (56) at (0, 3.1) {};
		\node [style=none] (57) at (3.1, 0) {};
		\node [style=none] (58) at (0, -3.1) {};
		\node [style=none] (59) at (-3.1, 0) {};
	\end{pgfonlayer}
	\begin{pgfonlayer}{edgelayer}
		\draw [style=circle, bend right=45, black, very thick] (56) to (59);
		\draw [style=circle, bend right=45, black, very thick] (57) to (56);
		\draw [style=circle, bend left=45, black, very thick] (57) to (58);
		\draw [style=circle, bend left=45, black, very thick] (58) to (59);
		\draw [style=circle, black, very thick] (7) to (14);
		\draw [style=circle, black, very thick] (14) to (42);
		\draw [style=circle, black, very thick] (42) to (28);
		\draw [style=circle, black, very thick] (7) to (28);
		\draw [style=circle] (15) to (55);
		\draw [style=circle] (54) to (16);
		\draw [style=circle] (17) to (53);
		\draw [style=circle] (52) to (18);
		\draw [style=circle] (19) to (51);
		\draw [style=circle] (50) to (20);
		\draw [style=circle] (21) to (49);
		\draw [style=circle] (22) to (48);
		\draw [style=circle] (47) to (23);
		\draw [style=circle] (24) to (46);
		\draw [style=circle] (45) to (25);
		\draw [style=circle] (26) to (44);
		\draw [style=circle] (43) to (27);
		\draw [style=circle] (13) to (41);
		\draw [style=circle] (40) to (12);
		\draw [style=circle] (39) to (11);
		\draw [style=circle] (38) to (10);
		\draw [style=circle] (37) to (9);
		\draw [style=circle] (36) to (8);
		\draw [style=circle] (35) to (0);
		\draw [style=circle] (34) to (1);
		\draw [style=circle] (33) to (2);
		\draw [style=circle] (32) to (3);
		\draw [style=circle] (31) to (4);
		\draw [style=circle] (30) to (5);
		\draw [style=circle] (29) to (6);
		\draw [style=circle] (6) to (15);
		\draw [style=circle] (5) to (16);
		\draw [style=circle] (4) to (17);
		\draw [style=circle] (18) to (3);
		\draw [style=circle] (2) to (19);
		\draw [style=circle] (1) to (20);
		\draw [style=circle] (0) to (21);
		\draw [style=circle] (8) to (22);
		\draw [style=circle] (9) to (23);
		\draw [style=circle] (10) to (24);
		\draw [style=circle] (11) to (25);
		\draw [style=circle] (12) to (26);
		\draw [style=circle] (27) to (13);
		\draw [style=circle] (14) to (28);
		\draw [style=circle] (29) to (55);
		\draw [style=circle] (30) to (54);
		\draw [style=circle] (31) to (53);
		\draw [style=circle] (52) to (32);
		\draw [style=circle] (33) to (51);
		\draw [style=circle] (34) to (50);
		\draw [style=circle] (35) to (49);
		\draw [style=circle] (36) to (48);
		\draw [style=circle] (37) to (47);
		\draw [style=circle] (38) to (46);
		\draw [style=circle] (39) to (45);
		\draw [style=circle] (40) to (44);
		\draw [style=circle] (41) to (43);
	\end{pgfonlayer}
\end{tikzpicture}
\caption{Example of mesh of $\Omega$.}
\label{domain_fictitious}
\end{center}
\end{figure}
Figure \ref{domain_fictitious} shows an example of configuration, the mesh do not fit with the interface $\Gamma$. Let
$$\Omega_i^*:=\left\{K\in\mathcal{T}_h~|~K\cap\Omega_i\neq \emptyset\right\},$$
We redefine the spaces defined in the previous section.
\begin{eqnarray*}
V_i&:=&\left\{v\in H^{1}\left(\Omega_i^*\right):v|_{\partial\Omega}=0\right\},\\
V_{i}^{kh}&:=&\left\{u_h\in V_i:u_h|_K\in\mathbb{P}_k\left(K\right)~~\forall K\in\mathcal{T}_h\right\},
\end{eqnarray*}
then $V_h^{k}:=V_{1}^{kh}\times V_{2}^{kh}$. We define the set of elements that intersects the boundary
$$G_h:=\left\{K\in\mathcal{T}_h~|~K\cap\Gamma\neq \emptyset\right\}.$$
For the sake of precision we make the following assumptions regarding the mesh $\mathcal{T}_h$ and the boundary $\Gamma$ :
\begin{itemize}
\item The boundary $\Gamma$ intersects each element boundary $\partial K$ exactly twice, and each (open) edge at most once for $K\in G_h$.
\item Let $\Gamma_{K,h}$ be the straight line segment connecting the points of intersection between $\Gamma$ and $\partial K$. We assume that $\Gamma_K$ is a function of length on $\Gamma_{K,h}$ ; in local coordinates
$$\Gamma_{K,h}=\left\{\left(\xi,\eta\right):0<\xi<\left|\Gamma_{K,h}\right|,\eta=0\right\}$$
and
$$\Gamma_{K}=\left\{\left(\xi,\eta\right):0<\xi<\left|\Gamma_{K,h}\right|,\eta=\delta\left(\xi\right)\right\}.$$
\item We assume that for all $K\in G_h$ there exists $K'\notin G_h$ and $K\cap K'\neq \emptyset$ and such that the measures of $K$ and $K'$ are comparable in the sense that there exists $c_q>1$ such that
$$c_q^{-1}\leq\frac{\text{meas}\left(K\right)}{\text{meas}\left(K'\right)}\leq c_q$$
and that the faces $F$ such that $K\cap F\neq \emptyset$ and $K'\cap F \neq \emptyset$ satisfy
$$\left|\text{meas}\left(K'\right)\right|\leq c_q \text{meas}\left(F\right)^2.$$
\item We assume that in a triangle $K$ intersected by the interface $\Gamma$, the scalar product between the normal $n_F$ of the face that does not intersects $\Gamma$ and the normal of the interface $n$ keeps the same sign in $K$.
\end{itemize}
We now extend the trace inequality for $v\in H^1\left(\Omega\right)$
\begin{equation}
\label{trace_fictitious}
\left\|v\right\|_{K\cap\Gamma}\lesssim \left(h_K^{-\frac12}\left\|v\right\|_{K}+h_K^{\frac12}\left\|\nabla v\right\|_{K}\right)~~~~~\forall K\in\mathcal{T}_h.
\end{equation}
The inequality \eqref{trace_fictitious} has been shown in \cite{Hansbo_2002_a}.
Let $\mathbb{E}_i$ be an $H^k$-extension on $\Omega_i^*$, $\mathbb{E}_i : H^{k}\left(\Omega_i\right)\rightarrow H^{k}\left(\Omega_i^*\right)$ such that $(\mathbb{E}_iw)|_{\Omega_i}=w$ and
\begin{equation}
\label{ineq_extension}
\left\|\mathbb{E}_i w\right\|_{H^k\left(\Omega_i^*\right)}\lesssim\left\|w\right\|_{H^k\left(\Omega_i\right)}~~~,~~~~~~~k\geq 0.
\end{equation}
Let $\pi_{h}^*$ be the standard nodal interpolant, we construct the interpolation operator $\pi_h$ such that $\pi_hw=(\pi_h^1w,\pi_h^2w)$ such that
\begin{equation}
\label{def_interpol}
\pi_h^i u^i:=\pi_h^*\mathbb{E}_i u^i
\end{equation}
We recall the interpolation estimates for $0\leq r \leq s \leq k+1$,
\begin{eqnarray}
\label{interpolation_estimate_lagrange}
\left\|u^i-\pi_h^*u^i\right\|_{H^r\left(K\right)}&\lesssim& h^{s-r} \left|u^i\right|_{H^{s}\left(K\right)}~~~~~~~~~\forall K\in \mathcal{T}_h,\\
\label{interpolation_estimate_lagrange_face}
\left\|u^i-\pi_h^*u^i\right\|_{H^r\left(F\right)}&\lesssim& h^{s-r-\frac12} \left|u^i\right|_{H^{s}\left(K\right)}~~~~~~\forall F\in\mathcal{F}.
\end{eqnarray}
Using the estimate (\ref{ineq_extension}) together with (\ref{interpolation_estimate_lagrange}) and (\ref{interpolation_estimate_lagrange_face})
\begin{eqnarray}
\label{interpolation_estimate}
\left\|u^i-\pi_h^iu^i\right\|_{H^r\left(\Omega_i^*\right)}&\lesssim& h^{s-r} \left|u^i\right|_{H^{s}\left(\Omega_i\right)},\\
\sum_{F\in \mathcal{F}}\left\|u^i-\pi_h^iu^i\right\|_{H^r\left(F\right)}&\lesssim& h^{s-r-\frac12} \left|u^i\right|_{H^{s}\left(\Omega_i\right)}.
\end{eqnarray}
The stability analysis for this case is similar to the fitted case treated previously, first we need to adapt the structure of patches to this new configuration. Let us split the set $G_h$ into $N_p$ smaller disjoint set of elements $G_j$ with $j=1,\dots,N_p$. Let $I_{G_j}$ be the set of nodes $\{x_n\}$ belonging to $G_j$. A generic node of the triangulation is designated by $x_n$, we define the sets of nodes $I_j^1$ and $I_j^2$ such that
$$I_j^1:=\left\{x_n\in I_{G_j}~|~x_n\in\Omega_1,~x_n\notin I_{G_i}~\forall i\neq j\right\},$$
$$I_j^2:=\left\{x_n\in I_{G_j}~|~x_n\in\Omega_2,~x_n\notin I_{G_i}~\forall i\neq j\right\},$$
now we define $P_j^1$ and $P_j^2$ for each $G_j$ such that
$$P_j^1:=G_j\cup \left\{K\in\mathcal{T}_h~|~I_j^1\in K\right\}~~~,~~~~~~P_j^2:=G_j\cup \left\{K\in\mathcal{T}_h~|~I_j^2\in K\right\}.$$
Each patch $P_j^i$ is constructed such that $I_j^i\neq\emptyset$ for $i=1,2$. Figure \ref{new_patches} shows an example of two patch $P_j^1$ and $P_j^2$ attached to a set of interface elements $G_j$. $\Gamma_j:=\Gamma \cap G_j$ is the part of the boundary included in the patches $P_j^1$ and $P_j^2$. For all $j$ and $i=1,2$, the patch $P_j^i$ has the following properties
\begin{equation}
\label{order_patch}
h\lesssim\mbox{meas}(\Gamma_j)\lesssim h ~~~~\text{and}~~~~h^2\lesssim\mbox{meas}(P_j^i)\lesssim h^2.
\end{equation}
The function $\phi_j$ attached to $P_j^1$ introduced in \eqref{phi_gamma} is modified such that 
\begin{equation*}
\phi_j\left(x_i\right)=
\left\{
\begin{array}{rllr}
0&\text{for}& x_i \notin I_j^1\\
1&\text{for}& x_i \in I_j^1,
\end{array}\right.
\end{equation*}
with $i=1,\dots,N_n$. $N_n$ is the number of nodes in the triangulation $\mathcal{T}_h$.
\begin{figure}[h!]
\begin{center}
\begin{tikzpicture}[scale=0.8]

%
%
%
%
%
%

	\begin{pgfonlayer}{nodelayer}
		\node [style=nfilled] (0) at (-2, -2) {};
		\node [style=nfilled] (1) at (-2, 0) {};
		\node [style=nfilled] (2) at (0, 2) {};
		\node [style=nfilled] (3) at (2, 2) {};
		\node [style=filled] (4) at (2, 0) {};
		\node [style=filled] (5) at (0, 0) {};
		\node [style=nfilled] (6) at (0, -2) {};
		\node [style=nfilled] (7) at (4, 4) {};
		\node [style=nfilled] (8) at (4, 2) {};
		\node [style=nfilled] (9) at (4, 0) {};
		\node [style=none] (10) at (-1.85, -2) {};
		\node [style=none] (11) at (4, 2.3) {};
		\node [style=nfilled] (12) at (2, -2) {};
	\end{pgfonlayer}
	\begin{pgfonlayer}{edgelayer}
		\draw [style=none] (0) to (5);
		\draw [style=none] (1) to (2);
		\draw [style=none] (5) to (3);
		\draw [style=none] (6) to (4);
		\draw [style=none] (4) to (8);
		\draw [style=none] (3) to (7);
		\draw [style=none] (2) to (3);
		\draw [style=none] (1) to (0);
		\draw [style=none] (0) to (6);
		\draw [style=none] (5) to (6);
		\draw [style=none] (1) to (5);
		\draw [style=none] (2) to (5);
		\draw [style=none] (5) to (4);
		\draw [style=none] (3) to (4);
		\draw [style=none] (3) to (8);
		\draw [style=none] (7) to (8);
		\draw [style=none] (8) to (9);
		\draw [style=none] (4) to (9);
		\draw [style=none, bend left, dashed, black, ultra thick] (10) to (11);
		\draw [style=none] (6) to (12);
		\draw [style=none] (4) to (12);
		\draw [style=none] (12) to (9);
	\end{pgfonlayer}
\end{tikzpicture}
\hspace{1cm}
\begin{tikzpicture}[scale=0.8]

%
%
%
%
%
%

	\begin{pgfonlayer}{nodelayer}
		\node [style=none] (0) at (-2, -2) {};
		\node [style=none] (1) at (-2, 0) {};
		\node [style=none] (2) at (0, 2) {};
		\node [style=none] (3) at (2, 2) {};
		\node [style=none] (4) at (2, 0) {};
		\node [style=none] (5) at (0, 0) {};
		\node [style=none] (6) at (0, -2) {};
		\node [style=none] (7) at (4, 4) {};
		\node [style=none] (8) at (4, 2) {};
		\node [style=none] (10) at (-1.85, -2) {};
		\node [style=none] (11) at (4, 2.3) {};
		\node [style=none] (13) at (-4, -2) {};
		\node [style=none] (14) at (-4, 0) {};
		\node [style=none] (15) at (-2, 2) {};
		\node [style=none] (16) at (0, 4) {};
		\node [style=none] (17) at (2, 4) {};
	\end{pgfonlayer}
	\begin{pgfonlayer}{edgelayer}
		\draw [style=none] (0) to (5);
		\draw [style=none] (1) to (2);
		\draw [style=none] (5) to (3);
		\draw [style=none] (4) to (8);
		\draw [style=none] (3) to (7);
		\draw [style=none] (2) to (3);
		\draw [style=none] (1) to (0);
		\draw [style=none] (0) to (6);
		\draw [style=none] (5) to (6);
		\draw [style=none] (1) to (5);
		\draw [style=none] (2) to (5);
		\draw [style=none] (5) to (4);
		\draw [style=none] (3) to (4);
		\draw [style=none] (3) to (8);
		\draw [style=none] (7) to (8);
		\draw [style=none, bend left, dashed, black, ultra thick] (10) to (11);
		\draw [style=none] (0) to (13);
		\draw [style=none] (14) to (13);
		\draw [style=none] (14) to (1);
		\draw [style=none] (13) to (1);
		\draw [style=none] (14) to (15);
		\draw [style=none] (1) to (15);
		\draw [style=none] (2) to (15);
		\draw [style=none] (2) to (16);
		\draw [style=none] (2) to (17);
		\draw [style=none] (16) to (17);
		\draw [style=none] (16) to (15);
		\draw [style=none] (17) to (3);
		\draw [style=none] (17) to (7);
	\end{pgfonlayer}
\end{tikzpicture}
\caption{Zoom of Figure \ref{domain_fictitious}. For a set $G_j$, left : example of $P_j^1$, the function $\phi_j$ is equal to 0 in the nonfilled nodes, 1 in the filled nodes ; right : example of $P_j^2$.}
\end{center}
\label{new_patches}
\end{figure}
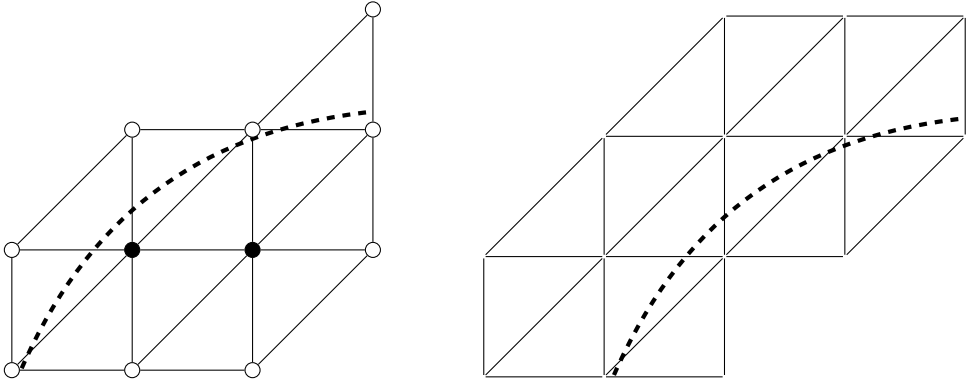
Similarly as the fitted case we define the function $v_\Gamma^1\in V_1^{kh}$ such that
\begin{equation}
\label{defvgamma2}
v_\Gamma^1=\alpha\sum_{j=1}^{N_P} v_j^1,
\end{equation}
with
\begin{equation}
\label{defvj2}
v_j^1=\nu_j\phi_j~~~,~~~~~~\nu_j\in\mathbb{R}.
\end{equation}
The function $v_j^1$ has the property
\begin{equation}
\label{vjprop2}
\text{meas}\left(\Gamma_j\right)^{-1}\int_{\Gamma_j}\nabla v_j^1\cdot n~\text{d}s=h^{-1}\overline{\llbracket u_h \rrbracket}^{\Gamma_j}.
\end{equation}
It is straightforward to remark that \eqref{stdapprox} and \eqref{poincare} still hold for this new configuration
\begin{equation}
\label{stdapprox2}
\left\|u_h^i-\overline{u_h^i}^{\Gamma_j}\right\|_{\Gamma_j}\lesssim h\left\|\nabla u_h^i\right\|_{\Gamma_j},
\end{equation}
\begin{equation}
\label{poincare2}
\left\|v_j^1\right\|_{P_j^1}\lesssim h\left\|\nabla v_j^1\right\|_{P_j^1}.
\end{equation}
The Lemma 1 of \cite{Boiveau_2015_c} combined with the regularity assumptions on the mesh made previously allows us to extend \eqref{keyineq_DD} to
\begin{equation}
\label{keyineq_DD2}
\left\|\nabla v_j^1\right\|_{P_j^1}^2\lesssim \left\|h^{-\frac12}\overline{\llbracket u_h\rrbracket}^{\Gamma_j}\right\|_{\Gamma_j}^2.
\end{equation}
Using the trace inequality \eqref{trace_fictitious}, Lemma \ref{lowbound_projection_DD} can be extended to
\begin{equation}
\sum_{j=1}^{N_p}\left\|\gamma^\frac12\overline{\llbracket u_h\rrbracket}^{\Gamma_j}\right\|_{\Gamma_j}\geq\sum_{j=1}^{N_p}\left\|\gamma^\frac12\llbracket u_h\rrbracket\right\|_{\Gamma_j}-C\omega_1\sum_{j=1}^{N_p}\left\|\mu_1^\frac12\nabla u_h^1\right\|_{P_j^1}-C\omega_2\sum_{j=1}^{N_p}\left\|\mu_2^\frac12\nabla u_h^2\right\|_{P_j^2}.
\label{lowbound_projection_DD2}
\end{equation}
We now have
$$\omega_1=\frac{\mu_2}{\mu_1+\mu_2}~~~,~~~~~~\omega_2=\frac{\mu_1}{\mu_1+\mu_2}~~~,~~~~~~\gamma=\frac{\mu_1\mu_2}{h(\mu_1+\mu_2)}.$$

\subsection{Finite element formulation}
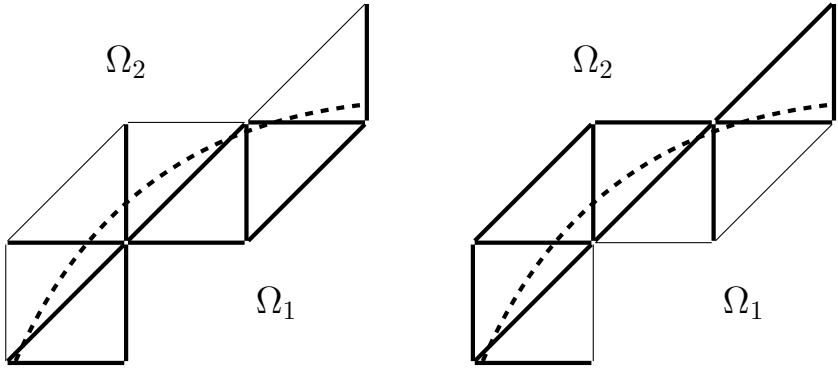
\begin{figure}[h!]
\begin{center}
\begin{tikzpicture}[scale=0.8]

%
%
%
%
%
%
\draw [black] (2.5,-1)node{\Large{$\Omega_1$}};
\draw [black] (0,3)node{\Large{$\Omega_2$}};
	\begin{pgfonlayer}{nodelayer}
		\node [style=none] (0) at (-2, -2) {};
		\node [style=none] (1) at (-2, 0) {};
		\node [style=none] (2) at (0, 2) {};
		\node [style=none] (3) at (2, 2) {};
		\node [style=none] (4) at (2, 0) {};
		\node [style=none] (5) at (0, 0) {};
		\node [style=none] (6) at (0, -2) {};
		\node [style=none] (7) at (4, 4) {};
		\node [style=none] (8) at (4, 2) {};
		\node [style=none] (9) at (4, 0) {};
		\node [style=none] (10) at (-1.85, -2) {};
		\node [style=none] (11) at (4, 2.3) {};
	\end{pgfonlayer}
	\begin{pgfonlayer}{edgelayer}
		\draw [style=none, black, ultra thick] (0) to (5);
		\draw [style=none] (1) to (2);
		\draw [style=none, black, ultra thick] (5) to (3);
		\draw [style=none, black, ultra thick] (4) to (8);
		\draw [style=none] (3) to (7);
		\draw [style=none] (2) to (3);
		\draw [style=none] (1) to (0);
		\draw [style=none, black, ultra thick] (0) to (6);
		\draw [style=none, black, ultra thick] (5) to (6);
		\draw [style=none, black, ultra thick] (1) to (5);
		\draw [style=none, black, ultra thick] (2) to (5);
		\draw [style=none,black, ultra thick] (5) to (4);
		\draw [style=none, black, ultra thick] (3) to (4);
		\draw [style=none, black, ultra thick] (3) to (8);
		\draw [style=none, black, ultra thick] (7) to (8);
		\draw [style=none, bend left, dashed, black, ultra thick] (10) to (11);
	\end{pgfonlayer}
\end{tikzpicture}
\hspace{1cm}
\begin{tikzpicture}[scale=0.8]

%
%
%
%
%
%
\draw [black] (2.5,-1)node{\Large{$\Omega_1$}};
\draw [black] (0,3)node{\Large{$\Omega_2$}};
	\begin{pgfonlayer}{nodelayer}
		\node [style=none] (0) at (-2, -2) {};
		\node [style=none] (1) at (-2, 0) {};
		\node [style=none] (2) at (0, 2) {};
		\node [style=none] (3) at (2, 2) {};
		\node [style=none] (4) at (2, 0) {};
		\node [style=none] (5) at (0, 0) {};
		\node [style=none] (6) at (0, -2) {};
		\node [style=none] (7) at (4, 4) {};
		\node [style=none] (8) at (4, 2) {};
		\node [style=none] (9) at (4, 0) {};
		\node [style=none] (10) at (-1.85, -2) {};
		\node [style=none] (11) at (4, 2.3) {};
	\end{pgfonlayer}
	\begin{pgfonlayer}{edgelayer}
		\draw [style=none, black, ultra thick] (0) to (5);
		\draw [style=none, ultra thick] (1) to (2);
		\draw [style=none, black, ultra thick] (5) to (3);
		\draw [style=none] (4) to (8);
		\draw [style=none, ultra thick] (3) to (7);
		\draw [style=none, ultra thick] (2) to (3);
		\draw [style=none, ultra thick] (1) to (0);
		\draw [style=none, black, ultra thick] (0) to (6);
		\draw [style=none] (5) to (6);
		\draw [style=none, black, ultra thick] (1) to (5);
		\draw [style=none, black, ultra thick] (2) to (5);
		\draw [style=none] (5) to (4);
		\draw [style=none, black, ultra thick] (3) to (4);
		\draw [style=none, black, ultra thick] (3) to (8);
		\draw [style=none, black, ultra thick] (7) to (8);
		\draw [style=none, bend left, dashed, black, ultra thick] (10) to (11);
	\end{pgfonlayer}
\end{tikzpicture}
\caption{Zoom of Figure \ref{domain_fictitious}, the dashed line is the boundary $\Gamma$, the bold facets are in the set $\mathcal{F}_G^i$, left $\mathcal{F}_G^1$, right $\mathcal{F}_G^2$.\label{fig:interior_facets}}
\end{center}
\end{figure}

To handle this problem we use the fictitious domain method on both subdomains $\Omega_1$ and $\Omega_2$. We assume $\mu_1\leq\mu_2$. Using the penalty free Nitsche's method, the finite element formulation for the problem \eqref{poisson_DD} is written as : find $u_h\in V_h^k$ such that
\begin{equation}
\label{formulation_fictitious}
A_h\left(u_h,v_h\right)+J_h\left(u_h,v_h\right)=L_h\left(v_h\right)~~~\forall v_h\in V_h^k,
\end{equation}
where
\begin{eqnarray*}
A_h\left(u_h,v_h\right)&=&\sum_{i=1}^2\left(\mu_i\nabla u_h^i,\nabla v_h^i\right)_{\Omega_i}-\left\langle\left\{\mu\nabla u_h \cdot n\right\},\left\llbracket v_h\right\rrbracket\right\rangle_{\Gamma}+\left\langle\left\{\mu\nabla v_h \cdot n\right\},\left\llbracket u_h\right\rrbracket\right\rangle_{\Gamma},\\
L_h\left(v_h\right)&=&\sum_{i=1}^2\left(f_i,v_h^i\right)_{\Omega_i}.
\end{eqnarray*}
The operator $J_h$ is the ghost penalty \cite{Burman_2010_a}, defined such that $J_h(u_h,v_h)=J_h^1\left(u_h^1,v_h^1\right)+J_h^2\left(u_h^2,v_h^2\right)$ with
\begin{eqnarray*}
J_h^1\left(u_h^1,v_h^1\right)&=&\gamma_g\sum_{F\in\mathcal{F}_G^1}\sum_{l=1}^k\left\langle \mu_1 h^{2l-1}\llbracket D_{n_F}^l u_h^1\rrbracket,\llbracket D_{n_F}^l v_h^1\rrbracket\right\rangle_F,\\
J_h^2\left(u_h^2,v_h^2\right)&=&\gamma_g\sum_{F\in\mathcal{F}_G^2}\sum_{l=1}^k\left\langle \mu_2 h^{2l-1}\llbracket D_{n_F}^l u_h^2\rrbracket,\llbracket D_{n_F}^l v_h^2\rrbracket\right\rangle_F.
\end{eqnarray*}
This penalization ensures the stability in the case of small cut elements that might cause a dramatic growth of the condition number. The sets $\mathcal{F}_G^i$ for $i=1,2$ are often called the interior facets
$$\mathcal{F}_G^1:=\left\{F\in G_h~|~F\cap\Omega^1\neq \emptyset\right\}~~~,~~~~~~\mathcal{F}_G^2:=\left\{F\in G_h~|~F\cap\Omega^2\neq \emptyset\right\}.$$
$D^l_{n_F}$ is the partial derivative of order $l$ in the direction $n_F$.
Figure \ref{fig:interior_facets} shows an example of the sets $\mathcal{F}_G^1$ and $\mathcal{F}_G^2$ for a small part of the boundary.
From  \cite{Massing_2014_a} we have the estimate
\begin{equation}
\label{ghostestimate}
\left\|\mu_i^\frac12\nabla v_h^i\right\|_{\Omega_i^*}^2\lesssim \left\|\mu_i^\frac12\nabla v_h^i\right\|_{\Omega_i}^2 + J_h^i(v_h^i,v_h^i)\lesssim \left\|\mu_i^\frac12\nabla v_h^i\right\|_{\Omega_i^*}^2.
\end{equation}
\subsection{Inf-sup stability}
We define the norm
$$\vvvert w \vvvert_{*}^2:=\sum_{i=1}^2\left\|\mu_1^{\frac12}\nabla w^i\right\|_{\Omega_i^*}^2+\left\|\gamma^\frac12\llbracket w \rrbracket\right\|_{\Gamma}$$
\begin{lemma}
For $u_h,v_h\in V_h^k$ with $v_h=u_h+v_\Gamma^1$, $v_\Gamma^1$ defined by equations \eqref{defvgamma2} and \eqref{defvj2} , there exists a positive constant $\beta_0$ such that the following inequality is true
$$\beta_0\vvvert u_h\vvvert_*^2\leq A_h(u_h,v_h)+J_h(u_h,v_h)$$
\end{lemma}
\proof
Using \eqref{ghostestimate} we can write
$$A_h(u_h,u_h)+J_h(u_h,u_h)=\sum_{i=1}^2(\mu_i\nabla u_h^i,\nabla u_h^i)_{\Omega_i}+J_h(u_h,u_h)\gtrsim\left\|\mu_1^\frac12\nabla u_h^1\right\|_{\Omega_1^*}^2+\left\|\mu_2^\frac12\nabla u_h^2\right\|_{\Omega_2^*}^2$$
\begin{equation*}
\begin{split}
(A_h+J_h)\left(u_h,v_j^1\right)=&~\left(\mu_1\nabla u_h^1,\nabla v_j^1\right)_{P_j^1\cap\Omega_1}-\left\langle\left\{\mu\nabla u_h\cdot n\right\}, v_j^1\right\rangle_{\Gamma_j}\\&+\omega_1\left\langle\mu_1 \nabla v_j^1\cdot n,\llbracket u_h\rrbracket\right\rangle_{\Gamma_j}+J_h\left(u_h,v_j^1\right).
\end{split}
\end{equation*}
Using the proof of Lemma \ref{lowerbound_poisson_DD} with \eqref{trace_fictitious} \eqref{keyineq_DD2}
\begin{eqnarray*}
\left(\mu_1\nabla u_h^1,\alpha\nabla v_j^1\right)_{P_j^1\cap\Omega_1}+\alpha J_h\left(u_h,v_j^1\right)
&\geq&-\left\|\mu_1^{\frac12}\nabla u_h^1\right\|_{P_j^1\cap\Omega_1}\alpha\mu_1^{\frac12}\left\|\nabla v_j^1\right\|_{P_j^1\cap\Omega_1}\\&&- J_h\left(u_h,u_h\right)^\frac12\alpha J_h\left(v_j^1,v_j^1\right)^\frac12\\
&\gtrsim&-\epsilon\left\|\mu_1^{\frac12}\nabla u_h^1\right\|_{P_j^1}^2-\frac{C\alpha^2\mu_1}{4\epsilon}\left\|\nabla v_j^1\right\|_{P_j^1}^2\\
&\gtrsim&-\epsilon\left\|\mu_1^{\frac12}\nabla u_h^1\right\|_{P_j^1}^2-\frac{C\alpha^2}{4\epsilon}\left(1+\frac{\mu_1}{\mu_2}\right)\left\|\gamma^{\frac12}\overline{\llbracket u_h\rrbracket}^{\Gamma_j}\right\|_{\Gamma_j}^2,
\end{eqnarray*}
Using the trace and inverse inequalities, \eqref{poincare2} and \eqref{keyineq_DD2} we have
\begin{equation*}
\left\langle\left\{\mu\nabla u_h\cdot n\right\}, \alpha v_j^1\right\rangle_{\Gamma_j}
\leq\frac{C\alpha^2}{2\epsilon} \left\|\gamma^\frac12\overline{\llbracket u_h\rrbracket}^{\Gamma_j}\right\|_{\Gamma_j}^2+\epsilon\omega_1\left\|\mu_1^\frac12\nabla u_h^1\right\|_{P_j^1}^2+\epsilon\omega_2\left\|\mu_2^\frac12\nabla u_h^2\right\|_{P_j^2}^2.
\end{equation*}
Using \eqref{vjprop2}, \eqref{keyineq_DD2}, \eqref{stdapprox2} the trace and the inverse inequalities we obtain
\begin{equation*}
\begin{split}
\alpha\omega_1\left\langle\mu_1\nabla v_j^1\cdot n,\llbracket u_h\rrbracket\right\rangle_{\Gamma_j}
\geq-&\epsilon\omega_1\left\|\mu_1^\frac12\nabla u_h^1 \right\|_{P_j^1}^2-\epsilon\omega_2\left\|\mu_2^\frac12\nabla u_h^2 \right\|_{P_j^2}^2\\&+\alpha\left(1-\frac{C\alpha}{2\epsilon}\right) \left\|\gamma^\frac12\overline{\llbracket u_h\rrbracket}^{\Gamma_j}\right\|_{\Gamma_j}^2.
\end{split}
\end{equation*}
We now have the lowerbound
\begin{equation*}
\begin{split}
A_h\left(u_h,v_h\right)\geq
\left\|\mu_1^{\frac12}\nabla u_h^1\right\|_{\Omega_1^*\backslash P^1}^2+&\left\|\mu_2^{\frac12}\nabla u_h^2\right\|_{\Omega_2^*\backslash P^2}^2+\sum_{j=1}^{N_p}\Bigg[\left(C_a-\omega_1CC_c\right)\left\|\mu_1^\frac12\nabla u_h^1\right\|_{P_j^1}^2\\&+\left(C_b-\omega_2CC_c\right)\left\|\mu_2^\frac12\nabla u_h^2\right\|_{P_j^2}^2+C_c\left\|\gamma^\frac12\llbracket u_h\rrbracket\right\|_{\Gamma_j}^2\Bigg],
\end{split}
\end{equation*}
with the constants
\begin{eqnarray*}
C_a&=&1-\epsilon\left(2\omega_1+1\right),\\
C_b&=&1-2\epsilon\omega_2,\\
C_c&=&\alpha\left(1-\alpha \frac{C}{4\epsilon}\left(5+\frac{\mu_1}{\mu_2}\right)\right).
\end{eqnarray*}
All terms are positive for $\epsilon=\text{min}\left[\frac{1}{2\left(2\omega_1+1\right)},\frac{1}{4\omega_2}\right]$ and $\alpha=\text{min}\left[\frac{4\epsilon\omega_1}{C\left(4\omega_1+1\right)},\frac{1}{2C}\right]$.
 \endproof
\begin{theorem}
\label{infsup_poisson_DD2}
There exists a positive constant $\beta>0$ such that for all functions $u_h\in V_h^k$ the following inequality holds
\begin{equation*}
\beta\left\vvvert u_h\right\vvvert_*\leq\underset{v_h\in V_h^k}{\textup{sup}} \frac{A_h\left(u_h,v_h \right)+J_h\left(u_h,v_h \right)}{\left\vvvert v_h\right\vvvert_*}.
\end{equation*}
\end{theorem}
\proof
Same as proof of Theorem \ref{infsup_poisson_DD}.
 \endproof
\subsection{A priori error estimate}
We have the following consistency relation
\begin{lemma}
\label{galerkin_poisson_DD2}
If $u\in H^{2}\left(\Omega_1\right)\times H^{2}\left(\Omega_2\right)$ is the solution of (\ref{poisson_DD}) and $u_h\in V_h^k$ the solution of (\ref{formulation_fictitious}) the following property holds
$$A_h\left(u-u_h,v_h\right)-J_h(u_h,v_h)=0~~,~~~~\forall v_h\in V_h^k.$$
\end{lemma}
Let us introduce the norm
$$\left\| w\right\|_*= \left\vvvert w\right\vvvert_*+\left\|\mu_1^\frac12 h^{\frac12}\nabla w^1\cdot n \right\|_{\Gamma}+\left\| \mu_2^\frac12 h^{\frac12}\nabla w^2\cdot n \right\|_{\Gamma}.$$
\begin{lemma}
\label{triplestar_poisson_DD2}
Let $w\in H^{2}\left(\Omega_1\right)\times H^{2}\left(\Omega_2\right)+V_h^k$ and $v_h\in V_h^k$, there exists a positive constant $M$ such that the bilinear form $A_h\left(\cdot,\cdot\right)$ has the property
$$A_h\left(w,v_h\right)\leq M\left\|w\right\|_* \left\vvvert v_h\right\vvvert_*.$$
\end{lemma}
\proof
Same as Lemma \ref{triplestar_poisson_DD} with $h_1=h_2=h$.
 \endproof
\begin{lemma}
\label{boundJ}
Let $u\in H^{2}\left(\Omega_1\right)\times H^{2}\left(\Omega_2\right)$ and $\pi_h$ the interpolation operator defined by (\ref{def_interpol}). For all $v_h\in V_h^k$ the following inequality holds
$$J_h\left(\pi_hu, v_h\right)\leq C_J h^k \left\|u\right\|_{H^{k+1}\left(\Omega\right)}\left\vvvert v_h \right\vvvert_*,$$
with ${\rm O}\left(C_J\right)={\rm O}\left(1\right)$.
\end{lemma}
\proof
The continuity of $u$ gives $J_h\left(\mathbb{E}u, v_h\right)=0$. Using the definition of the interpolant $\pi_h$ (\ref{def_interpol}), the inverse inequality (\ref{inverse}), the interpolation estimate (\ref{interpolation_estimate}) and the continuity of $\mathbb{E}$
$$J_h\left(\pi_hu, v_h\right)=J_h\left(\pi_h^*\mathbb{E}u-\mathbb{E}u, v_h\right)\lesssim h^k \left|\mathbb{E}u\right|_{H^{k+1}\left(\Omega_h\right)}\left\|\nabla v_h\right\|_{\Omega_h}\lesssim h^k\left\|u\right\|_{H^{k+1}\left(\Omega\right)}\left\vvvert v_h\right\vvvert_*.$$
 \endproof
\begin{proposition}
\label{bounderror_fictitious}
If $u\in H^{k+1}\left(\Omega_1\right)\times H^{k+1}\left(\Omega_2\right)$ is the solution of (\ref{poisson_DD}) and $u_h\in V_h^k$ the solution of (\ref{formulation_fictitious}), then there holds
$$ \left\vvvert u-u_h\right\vvvert \leq C_{f\mu} h^k\left|u\right|_{H^{k+1}\left(\Omega\right)}.$$
with $C_{f\mu}$ is a positive constant that depends on $\mu$ and the mesh geometry.
\end{proposition}
\proof
Using the orthogonality relation of Lemma \ref{galerkin_poisson_DD2} we can write
$$A_h\left(u_h-\pi_h u,v_h\right)+J_h(u_h-\pi_h u,v_h)=A_h\left(u-\pi_h u,v_h\right)-J_h(\pi_h u,v_h),$$
applying this property with Theorem \ref{infsup_poisson_DD2} and Lemma \ref{triplestar_poisson_DD2} we get
\begin{equation*}
\beta\left\vvvert u_h-\pi_h u\right\vvvert_*\leq \frac{A_h\left(u-\pi_h u,v_h\right)-J_h(\pi_h u,v_h)}{\left\vvvert v_h\right\vvvert_*}\leq M\left\| u-\pi_h u\right\|_*-\frac{J_h(\pi_h u,v_h)}{\left\vvvert v_h\right\vvvert_*}.
\end{equation*}
Using \eqref{interpolation_estimate_lagrange} we can write
\begin{equation*}
\left\vvvert u-\pi_hu\right\vvvert\lesssim\left\|u-\pi_hu\right\|_*
\lesssim \mu_1^\frac12 h^{k}\left|u^1\right|_{H^{k+1}\left(\Omega_1\right)}+\mu_2^\frac12 h^{k}\left|u^2\right|_{H^{k+1}\left(\Omega_2\right)}.
\end{equation*}
It is straightforward to observe that for any $v_h\in V_h^k$ we have $\vvvert v_h \vvvert \lesssim \vvvert v_h \vvvert_*$, using this result and the triangle inequality we have
\begin{eqnarray*}
\left\vvvert u-u_h \right\vvvert
&\leq& \left\vvvert u-\pi_h u \right\vvvert + \left\vvvert u_h-\pi_h u \right\vvvert_*\\
&\leq& \left\vvvert u-\pi_h u \right\vvvert + \frac1\beta\left(M\left\| u-\pi_h u\right\|_*-\frac{J_h(\pi_h u,v_h)}{\left\vvvert v_h\right\vvvert_*}\right).
\end{eqnarray*}
Applying Lemma \ref{boundJ}
$$ \left\vvvert u-u_h\right\vvvert \lesssim \left(1+\frac{M+1}{\beta}\right) h^k\left\|u\right\|_{H^{k+1}\left(\Omega\right)}. \endproof$$
\begin{lemma}
\label{l2cvg_fictitious}
Let $u\in H^{k+1}\left(\Omega_1\right)\times H^{k+1}\left(\Omega_2\right)$ be the solution of (\ref{poisson_DD}) and $u_h$ the solution of (\ref{formulation_fictitious}), then
$$\left\|u-u_h\right\|_{\Omega}\leq C_{f\mu}' h^{k+\frac12}\left|u\right|_{H^{k+1}\left(\Omega\right)}$$
with $C_{f\mu}'$ is a positive constant that depends on $\mu$ and the mesh geometry.
\end{lemma}
\proof
Same proof as Lemma \ref{l2cvg_poisson_DD} considering the new unfitted framework.
 \endproof

\section{Numerical verifications}
In this section we verify numerically the convergences proven theoretically. For each case studied the domain considered is the unit square separated in two subdomains as it is shown in Figure \ref{domain}. We use a manufactured solution in order to test the precision and determine the slopes of convergence. The manufactured solution that has been considered in this case
$$u=\text{exp}(xy)\text{sin}(\pi x)\text{sin}(\pi y).$$
For fitted and unfitted domain decomposition we consider $\mu_1=1$ and we test a range of value for $\mu_2$.

\subsection{Fitted domain decomposition}
The package FreeFem++ \cite{Hecht_2012_a} is used to implement this case. We choose different values of the ratio $h_1/h_2$ and observe the $L^2$ and $H^1$-error for each configuration.
\begin{figure}[H]
\begin{center}
\includegraphics[scale=0.27]{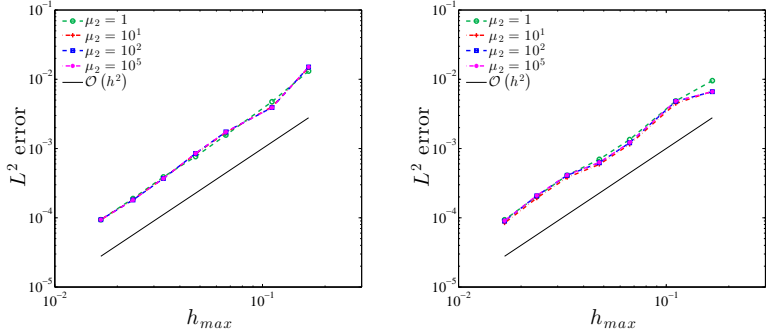}
\caption{Piecewise affine approximation, $\mu_1=1$, left $\frac{h_1}{h_2}=1$, right $\frac{h_1}{h_2}=\frac{3}{5}$.}
\label{poisson_L2_plot_1}
\end{center}
\end{figure}
\begin{figure}[H]
\begin{center}
\includegraphics[scale=0.27]{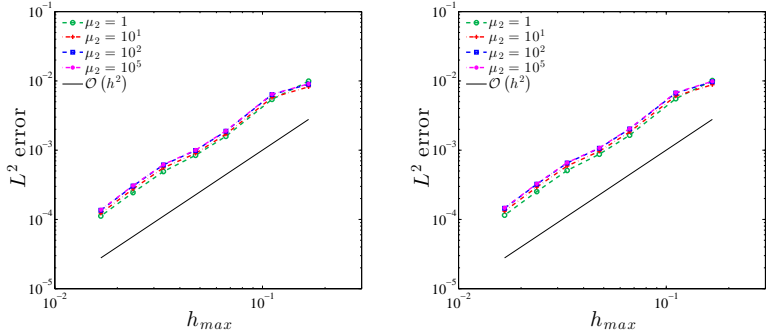}
\caption{Piecewise affine approximation, $\mu_1=1$, left $\frac{h_1}{h_2}=\frac{3}{10}$, right $\frac{h_1}{h_2}=\frac{1}{5}$.}
\label{poisson_L2_plot_2}
\end{center}
\end{figure}
Figure \ref{poisson_L2_plot_1} and \ref{poisson_L2_plot_2} shows a convergence of order $\mathcal{O}(h^2)$ for the $L^2$-error, this is a super convergence of order $\mathcal{O}(h^{\frac12})$ compared to the theoretical result. This super convergence has been observed for linear elasticity with the penalty free Nitsche's method in \cite{Boiveau_2015_a}. Comparing all four graphs we observe that as the ratio $h_1/h_2$ becomes smaller the constant $C_\mu'$ becomes slightly bigger when $\mu_2$ grows.

\begin{figure}[H]
\begin{center}
\includegraphics[scale=0.27]{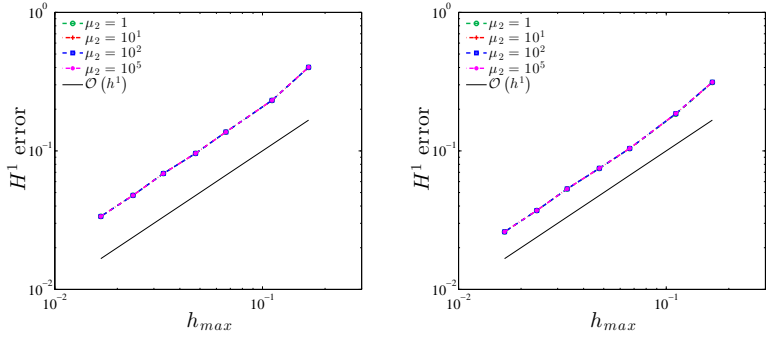}
\caption{Piecewise affine approximation, $\mu_1=1$, left $\frac{h_1}{h_2}=1$, right $\frac{h_1}{h_2}=\frac{3}{5}$.}
\label{poisson_H1_plot_1}
\end{center}
\end{figure}
\begin{figure}[H]
\begin{center}
\includegraphics[scale=0.27]{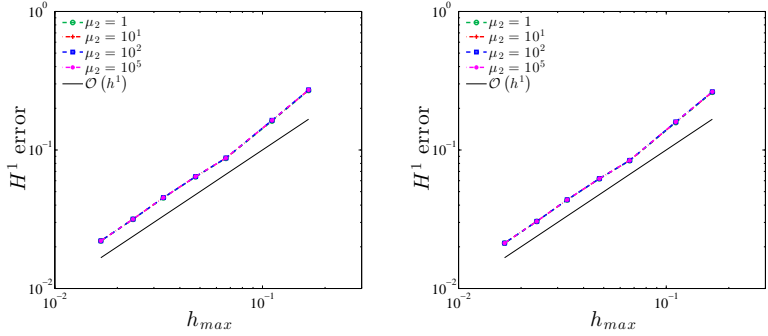}
\caption{Piecewise affine approximation, $\mu_1=1$, left $\frac{h_1}{h_2}=\frac{3}{10}$, right $\frac{h_1}{h_2}=\frac{1}{5}$.}
\label{poisson_H1_plot_2}
\end{center}
\end{figure}
We note that in all the cases the slope of convergence that has been proven theoretically is observed as it is shown in Figures \ref{poisson_H1_plot_1} and \ref{poisson_H1_plot_2}. In fact the affine approximation considered gives slopes of convergence of order $\mathcal{O}\left(h\right)$ which is what has been shown theoretically. For $h_1/h_2=1$ the meshsize are the same on both sides of $\Gamma$, in this case the influence of $\mu_2$ is negligible, the error remains the same for every value of $\mu_2$ considered. By considering the ratio $h_1/h_2$ smaller, the nonconformity of the meshes on both size of gamma gets bigger but it has a very small impact on the error for the three nonconforming cases considered.
\subsection{Unitted domain decomposition}
The package FEniCS \cite{Logg_2012_a} and the library CutFEM \cite{Burman_2014_g} have been used for these computations.

\begin{figure}[H]
\begin{center}
\includegraphics[scale=0.27]{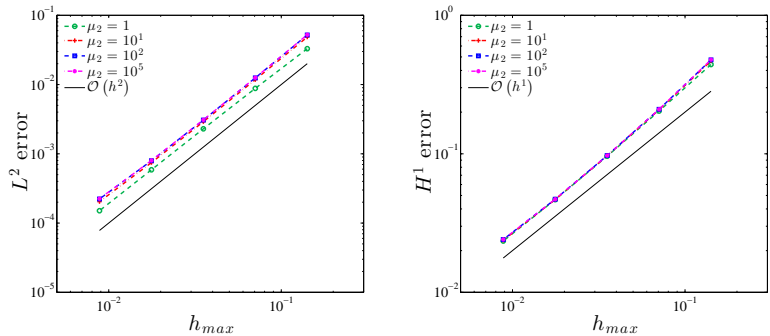}
\caption{Piecewise affine approximation, unfitted case, $\gamma_g=0.001$.}
\label{poisson_plot_unfitted}
\end{center}
\end{figure}

The same super convergence is observed as in the fitted case for the $L^2$-error. The case $\mu_2=1$ seems to have a constant $C_{f\mu}'$ slightly smaller than the other cases. The $H^1$-error shows the same convergence as shown in the theory, once again, the difference between each case is negligible.

	\bibliographystyle{siam.bst}
	\bibliography{Bibliography}

\end{document}